\documentclass[11pt]{article}

\usepackage{color}
\usepackage{latexsym}
\usepackage{dsfont}
\usepackage{amssymb}
\usepackage{comment}
\usepackage{graphicx}
\usepackage{amsmath, amsfonts,amssymb,theorem,euscript,array,enumerate,amsfonts,mathrsfs}
\usepackage{hyperref}
\usepackage{appendix}

\usepackage{stmaryrd} 
\usepackage{mathabx}

\usepackage{mathtools}
\mathtoolsset{showonlyrefs}

\usepackage{comment}

\usepackage{tikz}
\usetikzlibrary{fit,matrix,chains,positioning,decorations.pathreplacing,arrows}

\usepackage{mathtools}
\mathtoolsset{showonlyrefs}

\usepackage[a4paper,margin=1in,heightrounded]{geometry}

\usepackage{algpseudocode}
\usepackage{algorithm}

\usepackage[normalem]{ulem}

\newcommand{\independent}{\protect\mathpalette{\protect\independenT}{\perp}}
\def\independenT#1#2{\mathrel{\rlap{$#1#2$}\mkern2mu{#1#2}}}
\newcommand{\ba}{{\bf a}}

\def \b1{\bf{1}}

\def \N{\mathbb{N}}
\def \R{\mathbb{R}}

\def \E{\mathbb{E}}

\def \P{\mathbb{P}}

\def \T{\mathbb{T}}

\def \beps{\boldsymbol{\epsilon}} 
\def \bpi{\boldsymbol{\pi}} 
 
\def \bp{\boldsymbol{p}} 
\def \bpi{\boldsymbol{\pi}} 
 
\def \bF{\boldsymbol{F}} 
\def \bff{\boldsymbol{f}} 
\def \bU{\boldsymbol{U}} 
\def \bz{\boldsymbol{z}} 
 
\def \bd{\boldsymbol{d}}

\def \bnu{\boldsymbol{\nu}} 
\def \beps{\boldsymbol{\eps}} 
 
\def \bmrA{\boldsymbol{\mrA}} 
\def \bu{\boldsymbol{u}} 
 
\def \balpha{\boldsymbol{\alpha}} 
\def \bxi{\boldsymbol{\xi}} 
\def \bAc{\boldsymbol{\Ac}} 
\def \bVc{\boldsymbol{\Vc}} 
\def \mfa{\mathfrak{a}} 
 
\def \mra{\mathrm{a}}

\def \d{\mathrm{d}}

\def\esssup_#1{\underset{#1}{\mathrm{ess\,sup\, }}}

\def\argmin_#1{\underset{#1}{\mathrm{argmin\, }}}
\def\argmax_#1{\underset{#1}{\mathrm{argmax\, }}}

\def \Ac{{\cal A}}

\def \Fc{{\cal F}}
\def \Gc{{\cal G}}

\def \Kc{{\cal K}}
\def \Lc{{\cal L}}
\def \Pc{{\cal P}}

\def \Mc{{\cal M}}

\def \Oc{{\cal O}}
\def \Sc{{\cal S}}
\def \Tc{{\cal T}}
\def \Uc{{\cal U}}
\def \Vc{{\cal V}}
\def \Wc{{\cal W}}
\def \Xc{{\cal X}}
\def \Yc{{\cal Y}}

\def \eps{\varepsilon}

\def \ep{\hbox{ }\hfill$\Box$}

\def\reff#1{{\rm(\ref{#1})}}

\def\boX{{\boldsymbol X}}
\def\bolx{{\boldsymbol x}}

\def\boAc{{\boldsymbol \Ac}}

\def\boa{{\boldsymbol a}}

\def\bX{{\bf X}}

\def\be{{\bf e}}
\def\bx{{\boldsymbol x}}
\def\by{{\boldsymbol y}}

\def\bp{{\bf p}}

\def \mfa{\mathfrak{a}} 
\def \mra{\mathrm{a}} 
\def \mrA{\mathrm{A}}

\def \d{\mathrm{d}}

\def \bomrd{{\boldsymbol \d}}

\def\beqs{\begin{eqnarray*}}
\def\enqs{\end{eqnarray*}}
\def\beq{\begin{eqnarray}}
\def\enq{\end{eqnarray}}

\newcommand{\bl}[1]{\textcolor{blue}{#1}}

\newtheorem{Theorem}{Theorem}[section]
\newtheorem{Definition}{Definition}[section]
\newtheorem{Proposition}{Proposition}[section]

\newtheorem{Lemma}{Lemma}[section]

\newtheorem{Remark}{Remark}[section]

\numberwithin{equation}{section}

\begin{document}

%\title{
%Ad pricing in a large population\\with mean-field discrete choices 
%}

\author{
M\'ed\'eric MOTTE
\footnote{LPSM, Universit\'e Paris Cit\'e  \sf medericmotte at gmail.com
%The author acknowledges support of the DIM MathInnov. 
}
\qquad\quad
Huy\^en PHAM
\footnote{LPSM, Universit\'e Paris Cit\'e, and CREST-ENSAE, \sf pham at lpsm.paris
The author acknowledges support of the ANR 18-IDEX-001.  This work was also partially supported by the Chair Finance \& Sustainable Development / the FiME Lab (Institut Europlace de Finance)
}
}

\title{
Quantitative propagation of chaos for mean field Markov decision process with  common noise 
}

\maketitle

\abstract{
We investigate propagation of chaos for mean field  Markov Decision Process  with common noise (CMKV-MDP), and when  the optimization is performed over randomized open-loop controls on infinite horizon. 
We first  state  a rate of convergence of order $M_N^\gamma$, where $M_N$ is the mean rate of convergence in Wasserstein distance of the empirical measure, and $\gamma$ $\in$ $(0,1]$ is an explicit constant,  
in the limit  of the value functions of $N$-agent control problem with asymme\-tric  open-loop controls, 
towards the value function of CMKV-MDP.  Furthermore, we  show  how to explicitly construct  $(\epsilon+\Oc(M_N^\gamma))$-optimal policies for the $N$-agent model 
from $\epsilon$-optimal policies for the CMKV-MDP.  Our approach relies on  sharp comparison between the Bellman 
operators 
%fixed point equations 
in the $N$-agent problem and the CMKV-MDP, and fine coupling of 
empirical measures. 
%Finally, we  provide examples of application of the propagation of chaos result, by approximately solving several toy models for  $N$-agent targeted advertising problems with social influence via the resolution of the associated CMKV-MDP. 
}

%\vspace{5mm}

%\noindent {\bf Keywords:}  

%\vspace{5mm}

%\noindent {\bf JEL Classification}: 
%C70, C61

%\vspace{5mm}

%\noindent {\bf MSC Classification}: 
%91B26, 90B60, 60G55, 

%\newpage

%\tableofcontents

\section{Introduction}

We consider  a  social planner problem with $N$ cooperative agents in a mean-field discrete time model with common noise over an infinite horizon.  
The controlled state process $\boX$ $=$ $(X^i)_{i\in\llbracket 1,N\rrbracket}$ of  the $N$-agent model  is given by the dynamical random system 
\begin{equation} \label{dynXiN}
\left\{
\begin{array}{rcl}
	X^{i}_0 & = & x_0^i,   \\
	X^{i}_{t+1} & = &  F(X^{i}_t, \alpha^{i}_t, \frac{1}{N}\sum_{j=1}^N\delta_{(X^{j}_t,\alpha^{j}_t)},\eps^i_{t+1}, \eps^0_{t+1}),  \quad 
	t\in\N. 
\end{array}
\right.
\end{equation}
Here,  $x_0^i$, $i$ $\in$ $\llbracket 1,N\rrbracket$, are the initial states valued in a compact Polish space $\Xc$ with metric $d$,  $(\eps^i_t)_{i\in \llbracket 1, N\rrbracket t\in\N^\star}$ is a family of mutually i.i.d. random variables on some probability space 
$(\Omega, \Fc,\P)$,  valued in some measurable space $E$, and representing idiosyncratic noises, while $(\eps^0_t)_{t\in\N^\star}$  is another family of i.i.d. random variables valued  in some measurable space $E^0$, and representing the common noise (independent of idiosyncratic noise).  The control $\alpha^i$ followed by agent $i$, is a process, valued in some compact Polish space $A$ with metric $d_A$,  and  adapted with respect  to the filtration $(\Fc^N_t)_{{t\in\N}}$ generated by $\beps$ $=$ $\big((\eps^i_t)_{i\in \llbracket 1, N\rrbracket},\eps_t^0\big)_{t\in\N^\star}$ 
%$(\eps^i_t)_{i\in \llbracket 1, N\rrbracket t\in\N^\star}$,  $(\eps^0_t)_{t\in\N^\star}$, 
and also completed with a family of mutually  i.i.d. uniform random variables $\bU$ $=$ $(U^i_t)_{i\in\llbracket 1, N\rrbracket, t\in \N}$ that are used for randomization of the controls.  The mean-field interaction between the agents is formalized via the state transition function $F$ by the dependence upon the empirical measure of both  state/action 
of all the other agents: here $F$ is a measurable function from $\Xc \times A \times \Pc(\Xc\times A)\times E\times E^0$ into $\Xc$, where $\Pc(\Xc\times A)$ is the space of probability measures on the product space $\Xc\times A$. 
%satisfying some assumptions to be precised later. 

The objective of the social planner is to maximize over the set $\boAc$ of  $A^N$-valued  $(\Fc_t^N)_{{t\in\N}}$-adapted processes $\balpha$ $=$ $(\alpha_t^i)_{i\in\llbracket 1,N\rrbracket,t\in\N}$ a criterion in the form
\beqs
V_N^{\balpha}(\bolx_0) & :=& \E \Big[ \frac{1}{N} \sum_{i=1}^N  \sum_{t=0}^\infty \beta^t f \big(X_t^i,\alpha_t^i, \frac{1}{N}\sum_{j=1}^N\delta_{(X^{j}_t,\alpha^{j}_t)} \big) \Big], 
\enqs
where we set $\bolx_0$ $=$ $(x_0^i)_{i\in\llbracket 1,N\rrbracket}$ $\in$ $\Xc^N$ for the initial state of the $N$ agent system. Here $\beta$ $\in$ $(0,1)$ is a discount factor, and $f$ is a bounded measurable real-valued  function on  
$\Xc \times A \times \Pc(\Xc\times A)$. 
%satisfying some assumptions to be precised later. 
The value function for this optimization problem is defined on $\Xc^N$ as 
\begin{align}  \label{defVN}
V_N(\bolx_0) & :=  \;  \sup_{\balpha \in \boAc}  V_N^{\balpha}(\bolx_0),
\end{align} 
and we notice that problem \eqref{dynXiN}-\eqref{defVN} is a standard Markov Decision Process (MDP) with state space $\Xc^N$, action space $A^N$,  and (randomized) open-loop controls, and is the mathematical framework for reinforcement learning with multiple agents in interaction.

\vspace{1mm}

Let us now formulate  the asymptotic mean-field problem when the number of agents $N$ goes to infinity.  This consists formally in replacing empirical distributions by theoretical ones in the dynamic system and gain functions. 
The controlled state process $X$  of the representative agent is given by
\begin{equation} \label{dynXMKV}
\left\{
\begin{array}{rcl}
	X_0 & = & \xi_0, \\
	X_{t+1} & = &  F(X^{\alpha}_t, \alpha_t, \P^0_{(X^{}_t,\alpha_t)},\eps_{t+1}, \eps^0_{t+1}),  \quad 	t\in\N, 
\end{array}
\right.
\end{equation}
where we have renamed the uniform random sequence $(U^1_t)_{t\in\N}$  and  the noise $(\eps^1_t)_{t\in\N}$ by $(U_t)_{t\in\N}$ and $(\eps_t)_{t\in\N}$, and  the initial state $\xi_0$ is a $\Gc$-measurable random variable, with $\Gc$ a $\sigma$-algebra  independent of $(U_t)_{t\in\N}$, $(\eps_t)_{t\in\N}$, 
$(\eps^0_t)_{t\in\N^\star}$, with distribution law $\mu_0$ $\in$ $\Pc(\Xc)$ (the set of probability measures on $\Xc$). The control process $\alpha$ is an $A$-valued process, adapted with respect to the filtration generated by $\Gc$, $(U_t)_{t\in\N}$, $(\eps_t)_{t\in\N}$, 
$(\eps^0_t)_{t\in\N^\star}$, denoted by $\alpha$ $\in$ $\Ac$. Here  $\P^0$ and $\E^0$ represent  the conditional probability and expectation knowing the common noise $\eps^0$, and then, given a random variable $Y$, we denote by $\P^0_Y$ or $\Lc^0(Y)$ its conditional law knowing $\eps^0$.  The McKean-Vlasov (or mean-field)  control  problem consists in maximizing over  randomized open-loop controls $\alpha$ in $\Ac_{}$ the gain functional
\beqs
V^\alpha(\xi_0) & :=&   \E\Big[ \sum_{t=0}^\infty \beta^tf\big(X_t, \alpha_t, \P^0_{(X^{}_t,\alpha_t)}\big) \Big]. 
\enqs
The value function to this optimization problem is defined on $\Pc(\Xc)$ by 
\begin{align} \label{defVMKV} 
V(\xi_0) & : = \; \sup_{\alpha \in \Ac} V^\alpha(\xi_0), 
\end{align} 
and we recall from \cite{motte2019meanfield} that $V$ depends on $\xi_0$ only through its distribution (invariance in law), and we denote by misuse of notation: $V(\mu_0)$ $=$ $V(\xi_0)$. 
Problem \eqref{dynXMKV}-\eqref{defVMKV} is called  mean-field Markov Decision Process with common noise, or conditional McKean-Vlasov Markov Decision Process (CMKV-MDP in short),  with the peculiarity compared to standard MDP 
coming from the dependence of the state transition on the conditional distribution of the state/action.  In view of  propagation of chaos for particle systems usually derived for mean-field diffusion process (see \cite{szi89}), it  is  expected that CMKV-MDP provides a mean-field approximation of the $N$-agent MDP model.  

While the literature on mean-field control in continuous time, in particular the optimal control of McKean-Vlasov equations, is quite  important, see the monograph \cite{cardelbook2} for an overview and related references, there are rather few papers devoted to the discrete time framework. 
One of the first works is \cite{gasgau11} which studies the convergence of large interacting population process to a simple mean-field model when the state space is finite. 
The paper \cite{phawei16} studies a discrete-time McKean-Vlasov control problem with feedback controls on finite horizon, and derive the corresponding dynamic programming equation which is explicitly solved in the linear quadratic case. In \cite{carlautan19},  the authors consider mean-field control on infinite horizon with common noise with a discussion about  connections between closed-loop and open-loop policies, and propose $Q$-learning algortithms.  Our companion paper \cite{motte2019meanfield} 
deals with open-loop control and highlights the role of randomized controls with respect to standard Markov Decision Process (MDP). The value function is  characterized as a fixed point Bellman equation defined on the space of probability measures, and existence of 
$\epsilon$-optimal randomized feedback controls is proved.  The recent paper \cite{bauerle2021mean} studies mean-field control with deterministic closed-loop policies through the lens of MDP theory, and discusses the existence of optimal policies for the limiting mean-field problem as well as for the $N$-agent problem.

\vspace{1mm}

\noindent {\bf Main contributions.}  In this paper, we establish a quantitative propagation of result  for the $N$-agent MDP  towards the CMKV-MDP.  Our  contributions are twofold:
\begin{enumerate}
    \item[1.] We show in Theorem \ref{theo-chaosvalue} an explicit rate of convergence of the value functions under some assumptions to be precised later: there exists some positive constant $C$ (depending on the data of the problem) such that for all $\bolx$ $=$ $(x^i)_{i\in\llbracket 1,N\rrbracket}$ $\in$ $\Xc^N$, 
\beqs
\Big| V_N(\bolx) -  V\big(\frac{1}{N} \sum_{i=1}^N \delta_{x^i}\big) \Big| & \leq & C M_N^\gamma, 
\enqs
where $M_N$ is the mean rate of convergence in Wasserstein distance of the empirical measure (see \cite{fougui15}), and $\gamma$  $\in$ $(0,1]$ is an explicit constant depending on $\beta$ and $F$.  
\item[2.]  We prove  that any $\epsilon$-optimal  randomized feedback policy for the CMKV-MDP (including the case $\epsilon$ $=$ $0$, i.e.,  optimal randomized feedback policy whose existence is shown)  yields either an approximate optimal feedback control or an approximate randomized feedback control for the $N$-agent MDP problem, in a constructive sense to be precised later with an explicit rate of convergence, see Theorems \ref{theo1-controlNMDP} and \ref{theo2-controlNMDP}. 
\end{enumerate}
While the first statement  for convergence of value function is important in theory, the second statement is particularly interesting in practice (but often less studied in the literature) since it means that if the McKean-Vlasov MDP is simpler to solve than the $N$-agent MDP (some examples and applications to targeted advertising are developed in the 
PhD thesis \cite{mottethesis21}), then one can compute an almost optimal randomized feedback policy for the McKean-Vlasov MDP, and then use it in the $N$-agent MDP: this will guaranty us to have an almost optimal control.

\vspace{1mm}

\noindent {\it Related literature.}  The convergence of the $N$-individual problem to the limiting mean-field control problem has been first rigorously proved in \cite{lac17} by tightness and martingale arguments for continuous-time controlled McKean-Vlasov equations. This result  has been extended in \cite{dje20} to the common noise case and when there is interaction via the joint distribution of the state and control.  The paper \cite{ganetal21}  proved by viscosity solutions method via the characterization of the Hamilton-Jacobi-Bellman equation  the convergence of the value function towards the $N$-agent problem to the value function of the mean-field control problem in the common noise case but without idiosyncratic noise.  Rate of convergence of order $1/N$ has been stated in \cite{gerphawar21} by Backward Stochastic Differential Equations techniques but under the strong condition that there exists a smooth solution to the Master Bellman equation. The recent paper \cite{caretal22} removed this regularity assumption on the value function, and obtained an algebraic rate of convergence of order $N^{-\gamma}$ for some constant $\gamma$ $\in$ $(0,1]$.  We mention also in the continuous-time framework the paper \cite{cec21} which derived a rate of convergence of order $N^{-1/2}$ when the state space is finite. 

%\red{mostly continuous time, our paper AAP for symmetric policies.
%Literature on CV.  
%We point out  the work \cite{lac17}, which is the first paper to rigorously connect mean-field control to large systems of controlled processes, 
%see also the recent paper  \cite{foretal18}.  mention that convergence is usually stated for $\E[V_N(\xi^1,\ldots,\xi^N]$, and no rate of CV for the approximate policies. Here constructive approach and explicit rate. 
%}

%\vspace{1mm}

The convergence of the value function  in the  $N$-agent problem in a discrete-time mean field framework has been studied in our companion paper \cite{motte2019meanfield}. However, it was assumed there that each agent used the same open-loop policy, applied to her own idiosyncratic noise and the common noise. In particular,  
agent's controls cannot depend upon other agent's idiosyncratic noises, and they have symmetric (or exchangeable) behaviours. This restriction was crucial for using  propagation of chaos argument relying on a pathwise comparison between the state and control 
processes in the $N$-individual model and the McKean-Vlasov MDPs.

In this paper, we consider that the control of each agent can also depend upon the idiosyncratic noises of all the population, and that they can do so in a completely asymmetric way (i.e. each agent can use a different open-loop policy). This additional flexibility and generality in the definition of controls prevents us 
from coupling controls between the $N$-agent and the McKean-Vlasov MDPs in a one-to-one fashion as in \cite{motte2019meanfield}. In order to overcome this difficulty, we adopt quite different arguments by coupling the Bellman operators instead of the state/control process of the $N$-agent and CMKV MDPs. More precisely, the 
strategy of the proof is the following: 

\noindent {\it Idea of the proof.} 
\begin{enumerate}
\item[(i)] We first derive the Bellman equation for the $N$-agent MDP, with arguments similar to \cite{motte2019meanfield}, i.e. we prove that $\Tc_N V_N = V_N$, where $\Tc_N$ is the operator defined by
\begin{align}
\Tc_N W(\bx) & := \; \sup_{\ba \in A^N} \T^{\ba}_N W (\bx),  \quad \bx \in \Xc^N,  \end{align}
with
\begin{align}
\T^{\ba}_N W (\bx) & := \; \frac{1}{N}\sum_{i=1}^N f(x^i,a^i, \frac{1}{N} \sum_{j=1}^N \delta_{(x^j,a^j)} ) + \beta \E \big[ W\big(  (F(x^i,a^i,\frac{1}{N} \sum_{j=1}^N \delta_{(x^j,a^j)}, \eps^i_1,\eps^0_1)_{i\in\llbracket 1, N\rrbracket}\big) \big],  
\end{align}
for $\bx=(x^i)_{i\in \llbracket 1,N\rrbracket} \in \Xc^N$, $\ba = (a^i)_{i\in \llbracket 1,N\rrbracket} \in A^N$.  
This property is obtained by seeing the $N$-agent MDP as a standard MDP on $\Xc^N$ with actions space $A^N$.
\item[(ii)]  Then, we observe that the operators $\T^\mra$ of the McKean-Vlasov MDP, derived in \cite{motte2019meanfield}, are, formally, the limits of $\T^{\ba}_N$ when $N\rightarrow \infty$, for $\mra\in L^0(\Xc\times [0,1], A)$ and $\ba\in A^N$ well coupled. Inspired by this formal observation, we ``compare'' $\T^{\ba}_N$ to $\T^\mra$ and prove that they are indeed ``close'' in some sense, for $N$ large. A key point is that $\T^{\ba}_N$ is defined on $L^\infty_m(\Xc^N)$ (the set of bounded measurable functions on $\Xc^N$, valued in $\R$) while $\T^a$ is defined on $L^\infty_m(\Pc(\Xc))$ (the set of bounded measurable functions on $\Pc(\Xc)$, valued in $\R$). To compare both type of objects, we introduce a canonical way to associate to a function $W\in L^\infty_m(\Pc(\Xc))$ the function $\widecheck{W}\in L^\infty_m(\Xc^N)$ by setting $\widecheck{W}(\bx)=W\big(\frac{1}{N} \sum_{i=1}^N \delta_{x^i} \big)$.
\item[(iii)] Once the proximity between $\T^{\ba}_N$ and $\T^{a}$ is established in a general sense, we prove the proximity of the value functions $V_N$ and $V$ by seeing them as the unique fixed points of the Bellman operators  
$\Tc_N=\sup_{\ba\in A^N}\T^{\ba}_N$ and $\Tc=\sup_{\mra\in L^0(\Xc\times [0,1], A)}\T^{\mra}$, following the intuition  that if two contracting operators are close, their unique fixed points should also be close.
\item[(iv)] Finally, we provide two procedures to build $\Oc(\epsilon+M_N^\gamma)$-optimal policies for the $N$-agent MDP from an $\epsilon$-optimal stationary randomized feedback policy for the McKean-Vlasov MDP. The idea is to view, for each MDP, any $\epsilon$-optimal policy as a policy satisfying the verification theorem, which is a property only linked to the Bellman operator, again following the intuition  that if two Bellman operators are close, the policies satisfying their verification results should also be close.
\end{enumerate}

\vspace{1mm}

\noindent {\bf  Outline of the paper.}   The rest of the paper is organized  as follows. We state the assumptions and the main results in Section \ref{secmain}, while Section \ref{sec-chaos-N} is devoted to their proofs.  Finally, we give in 
Appendix \ref{sec-optcon} the proof of existence for optimal randomized feedback policy, and put in Appendix \ref{sec-Bell-N} some results about the Bellman operator for the $N$-agent MDP problem that are 
needed in the proof of our convergence results.

 \section{Main results} \label{secmain}

\subsection{Notations and assumptions} 
 
The  product space $\Xc\times A$ is  equipped with the metric $\bomrd((x,a),(x',a'))$ $=$ $d(x,x')$ $+$ $d_A(a,a')$, $x,x'$ $\in$  $\Xc$, $a,a'$ $\in$ $A$. 
%and the associated space of probability measure $\Pc(\Xc\times A)$ is equipped with its Wasserstein distances $\bW$. 
Likewise, we shall endow $\Xc^N$ with the metric $\bd_N(\bx,\bx')=\frac{1}{N}\sum_{i=1}^N d(x^i,x'^i)$ for $\bx=(x^i)_{i\in\llbracket 1,N\rrbracket},\bx'=(x'^i)_{i\in\llbracket 1,N\rrbracket}\in\Xc^N$, 
$A^N$ with the metric $\bd_{A,N}(\ba,\ba')=\frac{1}{N}\sum_{i=1}^N d_A(a^i,a'^i)$ for $\ba=(a^i)_{i\in\llbracket 1,N\rrbracket},\ba'=(a'^i)_{i\in\llbracket 1,N\rrbracket}\in A^N$, and $(\Xc\times A)^N$ with the metric $\bomrd_{N}((\bx,\ba),(\bx',\ba'))=\frac{1}{N}\sum_{i=1}^N \bomrd((x^i,a^i),(x'^i,a'^i))$ for $\bx, \bx'\in\Xc$ and $\ba,\ba'\in A^N$.
When $(\Yc,d)$ is a compact metric space, the set $\Pc(\Yc)$ of probability measures on $\Yc$ is equipped with the Wasserstein distance
\begin{align}
\Wc_d(\mu,\mu') & := \inf\Big\{ \int_{\Yc^2} d(y,y') \boldsymbol{\mu}(\d y,\d y'):  \boldsymbol{\mu} \in \boldsymbol{\Pi}(\mu,\mu') \Big\},
%\Pc(\Yc\times\Yc) \mbox{ with marginals } \mu \mbox{ and } \mu' \Big\}.  
%\\ \Wc_A(\nu,\nu') &= \inf\Big\{ \int d_A(a,a') \boldsymbol{\nu}(da,da'):  \boldsymbol{\mu} \in \Pc(A\times A) \mbox{ with marginals } \nu \mbox{ and } \nu' \Big\}, 
\end{align} 
where $\boldsymbol{\Pi}(\mu,\mu')$ is the set of (coupling) probability measures on $\Yc\times\Yc$ with marginals $\mu$ and $\mu'$, and we recall the dual Kantorovich-Rubinstein representation 
\begin{align} \label{dualW}
\Wc_d(\mu,\mu') & = \; \sup_{\phi \in Lip_1} \int_\Yc \phi(y) (\mu-\mu')(\d y),
\end{align} 
where $Lip_1$ is the set of Lipschitz functions on $\Yc$ with Lipschitz constant bounded by $1$. 

Given $\bolx$ $=$ $(x^i)_{i\in\llbracket 1,N\rrbracket}$ $\in$ $\Xc^N$, and $\ba=(a^i)_{i\in\llbracket 1,N\rrbracket}$ $\in$ $A^N$, we denote by 
\beqs
\mu_{_N}[\bolx] \; : = \;  \frac{1}{N} \sum_{i=1}^N \delta_{x^i} \; \in \; \Pc(\Xc), & & \mu_{_N}[\bolx,\ba] \; : = \;  \frac{1}{N} \sum_{i=1}^N \delta_{(x^i,a^i)} \; \in \; \Pc(\Xc\times A),
\enqs
and we recall that 
\begin{align} \label{inegW} 
\Wc_{\bomrd}(\mu_{_N}[\bx,\ba], \mu_{_N}[\bx',\ba']) & \leq \; \bomrd_{N}((\bx,\ba),(\bx',\ba')).
\end{align}

Given a random variable $Y$ on $(\Omega,\Fc,\P)$,  we denote  by $\P_Y$ or $\Lc(Y)$ its distribution law.

 \vspace{3mm}

 We make the following standing  assumptions on the  state transition function $F$ and on the running reward function $f$.

\vspace{3mm}

\noindent $({\bf HF_{lip}})$ There exists $K_F$ $>$ $0$,  such that   for all $a,a'\in A$, $e^0$ $\in$ $E^0$, $x,x'\in\Xc$, $\mu,\mu'\in\Pc(\Xc\times A)$,
\begin{align*}
\E \big[ d\big(F(x,a,\mu,\eps^1_1, e^0), F(x',a',\mu',\eps^1_1, e^0)\big)\big] & \leq  \; K_F \big( \bomrd((x,a),(x',a')) + \Wc_{\bomrd}(\mu,\mu')  \big)). 
\end{align*}

\vspace{1mm} 

\noindent $({\bf Hf_{lip}})$ There exists $K_f$ $>$ $0$,  such that   for all  $x,x'\in\Xc$, $a,a'\in A$, $\mu,\mu'\in\Pc(\Xc\times A)$,
\begin{align*}
\vert f(x,a,\mu)- f(x',a',\mu')\vert  & \leq  \; K_f  \big( \bomrd((x,a),(x',a')) + \Wc_{\bomrd}(\mu,\mu')  \big). 
\end{align*}

 {\begin{Remark}
\rm{
We stress the importance of making the regularity assumptions for $F$ in {\em expectation} only. When $\Xc$ is finite, $F$ cannot be, strictly speaking, 
Lipschitz (or even continuous) unless it is constant w.r.t. its mean-field argument ($\mu$ and $\mu'$ in $({\bf HF_{lip}})$). However, $F$ can be Lipschitz {\em in expectation}, e.g. once integrated w.r.t. the idiosyncratic noise. 
%which is a very natural assumption, as we shall see in Section \ref{sec-examples-finite}.
} 
%For example, setting $\Xc=\{0,1\}$ (and identifying $\Pc(\Xc)$ to the set of Bernoulli parameters $[0,1]$), and setting the idiosyncratic noise to be uniform in $[0,1]=:E$, 
%the function $F:(p,e)\in\Pc(\Xc)\times E\mapsto {\bf 1}_{e<p}\in\Xc$ is obviously not Lipschitz, but $\E [ d( F(p,\eps_1), F(p',\eps_1)) ] \leq  \vert p-p'\vert$, so $F$ is Lipschitz in expectation w.r.t. the idiosyncratic noise.
%An simple example of such function is given in Remark \ref{lipdiscex}.
%\ep
\end{Remark}

\vspace{2mm}

Under Assumption $({\bf HF_{lip}})$, we define the constant 
\beqs
\gamma & := &   \min\big[ 1, \frac{\vert \ln \beta\vert}{\ln(2K_F)_+}\big] \; \in \; (0,1]. 
\enqs

\vspace{3mm}

In the sequel, we denote by $\Delta_{\Xc}$ (resp. $\Delta_A$ and $\Delta_{\Xc\times A}$) the diameter of the compact metric space $\Xc$ (resp. $A$ and  $\Xc\times A$), and  define 
\begin{align} \label{defMN}
M_N &:= \underset{\mu\in\Pc(\Xc\times A)}{\sup}\E [\Wc_{\bomrd}(\mu_N,\mu)],
\end{align}
 where $\mu_N$ is the empirical measure $\mu_N$ $=$ $\frac{1}{N}\sum_{n=1}^N\delta_{Y_n}$, $(Y_n)_{1\leq n\leq N}$ are i.i.d. random variables with law $\mu$.  It is know that $M_N\underset{N\rightarrow \infty}{\rightarrow}0$, and 
we recall  from \cite{fougui15}, and  \cite{boileg14}  some results about non asymptotic bounds for  the mean rate of convergence in Wasserstein distance of the empirical measure. 
\begin{itemize}
\item If $\Xc\times A\subset \R^d$ for some $d\in\N^\star$, then:
$M_N=\Oc(N^{-\frac{1}{2}})$	for $d=1$, $M_N=\Oc(N^{-\frac{1}{2}}\log (1+N))$	for $d=2$, and $M_N=\Oc(N^{-\frac{1}{d}})$	for $d\geq 3$.
\item If for all $\delta>0$, the smallest number of balls with radius $\delta$ covering the compact metric set $\Xc\times A$ with diameter $\Delta_{\Xc\times A}$  
is smaller than $\Oc\Big( \big(\frac{\Delta_{\Xc\times A}}{\delta}\big)^\theta\Big)$ for $\theta>2$, then $M_N= \Oc(N^{-1/\theta})$.
\end{itemize}

 \vspace{3mm}
 
 In the sequel $C$ will denote a generic constant that depends only on the data of the problem, namely $\Delta_\Xc$, $\Delta_{\Xc\times A}$, $\beta$, $K_F$ and $K_f$.

\subsection{Convergence of value functions}

Our first main result is to quantify the rate of convergence of the value function of the $N$-agent MDP towards the value function of the CMKV-MDP.

\begin{Theorem}\label{theo-chaosvalue}
%Under Assumptions  $({\bf HF_{lip}})$  and $({\bf Hf_{lip}})$, we set  $\gamma = \min\big[ 1, \frac{\vert \ln \beta\vert}{\ln(2K_F)_+}\big]$.  Then, 
There exists some  positive constant $C$ 
%(depending on $\Delta_{\Xc\times A}$, $\beta$, $K_F$, $K_f$  \red{to check}) 
such that for all $\bolx$ $=$ $(x^i)_{i\in\llbracket 1,N\rrbracket}$ $\in$ $\Xc^N$, 
\beqs
\Big| V_N(\bolx) -  V\big(  \mu_{_N}[\bolx] \big) \Big| & \leq & C M_N^\gamma. 
\enqs
\end{Theorem}

\subsection{Approximate optimal policies}

Our next results are to show how to obtain approximate optimal control for the $N$ agent MDP from $\eps$-optimal control for CKMV-MDP, and to quantify  the accuracy of this approximation.  

First, let us recall from \cite{motte2019meanfield} the construction of $\eps$-optimal control for CKMV-MDP.  The value function $V$ is characterized as the unique fixed  point  in $L_m^\infty(\Pc(\Xc))$,  the set of bounded measurable real-valued functions on $\Pc(\Xc)$, of the Bellman equation $V$ $=$ $\Tc V$, where $\Tc$ is the Bellman operator defined on $L_m^\infty(\Pc(\Xc))$ by 
\begin{align}
\Tc  W(\mu) & := \;  \sup_{\mra \in  L^0(\Xc\times [0,1]; A)} \T^\mra W(\mu),  \\
\mbox{ with } \quad \T^\mra W(\mu) & := \;  
\E \Big[ f(\xi,\mra(\xi,U),\Lc(\xi,\mra(\xi,U))) + \beta W\big(  \P^0_{F(\xi,\mra(\xi,U),\Lc(\xi,\mra(\xi,U)),\eps_1,\eps_1^0)}   \big) \Big], \label{defTa} 
\end{align} 
for any $(\xi,U)$ $\sim$  $\mu\otimes\Uc([0,1])$ (it is clear that the right-hand side in \eqref{defTa} does not depend on the choice of such $(\xi,U)$), where $L^0(\Xc\times [0,1]; A)$ is the set of measurable functions from $\Xc\times [0,1]$ into $A$. 
Then, for all $\epsilon$ $>$ $0$, there exists a randomized feedback policy $\mfa_\epsilon$, i.e. a measurable function from $\Pc(\Xc)\times \Xc\times [0,1]$ into $A$,  denoted by $\mfa_\epsilon$ $\in$ $L^0(\Pc(\Xc)\times \Xc\times [0,1];A)$, 
%denoted  $\mfa_\eps$ $\in$   $L^0(\Pc(\Xc)\times \Xc\times [0,1];A)$, 
such that for all $\mu$ $\in$ $\Pc(\Xc)$:
\beqs
V(\mu) - \epsilon & \leq & \T^{\mfa_\epsilon(\mu,.)} V(\mu),  
\enqs
and we say that $\mfa_\epsilon$ is an $\epsilon$-optimal randomized feedback policy for CMKV-MDP.  
By considering the randomized feedback control $\alpha^\epsilon$ $\in$ $\Ac$ defined by 
\begin{align} \label{ranfeedbackeps} 
\alpha_t^\epsilon &= \;  \mfa_\epsilon(\P_{X_t}^0,X_t,U_t), \quad t \in \N, 
\end{align} 
where $(U_t)_{t\in\N}$ is an i.i.d. sequence of random variables, $U_t$ $\sim$ $\Uc([0,1])$, independent of $\xi_0$ $\sim$ $\mu_0$, and $\beps$,  this yields an $O(\epsilon)$-optimal control for $V(\mu_0)$,  namely
\beqs
V(\mu_0)  - \frac{\epsilon}{1-\beta} & \leq &  V^{\alpha^\epsilon}(\xi_0). 
\enqs
Actually, we can even take $\epsilon$ $=$ $0$, i.e.,  get optimal randomized feedback control.  The proof for the existence of an optimal randomized feedback policy 
is inspired by the paper \cite{carlautan19}, which states the existence of an optimal policy in a closely related model,  
%follows similar arguments as in \cite{carlautan19}, 
and is reported in Appendix \ref{sec-optcon}.

\vspace{1mm}

We now provide two procedures to construct approximate optimal control for the $N$-agent MDP from an $\epsilon$-optimal randomized feedback policy for CMKV-MDP.  The first procedure gives a general approach for getting approximate feedback control for the $N$-agent MDP.

\begin{Theorem} \label{theo1-controlNMDP}
Let $\mfa_\epsilon$ be an $\epsilon$-optimal randomized feedback policy for CMKV-MDP.  Then, there exists a measurable function $\bpi^{\mfa_\epsilon,N}$ from $\Xc^N$ into $A^N$, called feedback policy for the $N$-agent MDP,  such that 
\begin{align}  \label{argmina}  
\bpi^{\mfa_\epsilon,N}(\bx)  & \in \;   \argmin_{\ba\in A^N}\Wc_{\bomrd}\big(\Lc \big( \xi_\bx,\mfa_\epsilon(\mu_{_N}[\bx],\xi_\bx,U) \big), \mu_{_N}[\bx,\ba] \big), \quad \bx \in \Xc^N, 
\end{align} 
with $(\xi_\bx,U)$ $\sim$ $\mu_{_N}[\bx]\otimes\Uc([0,1])$. This yields a feedback control $\balpha^{\epsilon,N}$ $\in$ $\boAc$ defined by 
\beqs
\balpha_t^{\epsilon,N} &=&  \bpi^{\mfa_\epsilon,N}(\boX_t), \quad t \in \N,   
\enqs
which is $O(\eps + M_N^\gamma)$-optimal control for $V_N(\bx_0)$, namely: 
%there exists some  positive constant $C$  (depending on $\Delta_{\Xc\times A}$, $\beta$, $K_F$, $K_f$) s.t. 
\beqs 
 V_N(\bx_0)  -  C \big[ \epsilon +  M_N^\gamma] & \leq & V_N^{\balpha^{\epsilon,N}} (\bx_0).
 \enqs
\end{Theorem}

\vspace{1mm}

Theorem  \ref{theo1-controlNMDP} provides a generic way to obtain a $\Oc(\epsilon+M_N^\gamma)$-optimal  feedback policy  for the $N$-agent MDP from an $\eps$-optimal  randomized feedback policy $\mfa_\eps$ for CMKV-MDP, simply by sending actions $\ba$ $=$ $(a^i)_{i\in\llbracket 1,N\rrbracket}$ 
to the population so that, once in state $\bx$,  the state-action pair $(\bx,\ba)$ is  empirically distributed as closely as possible to $\Lc\big(\xi_\bx,\mfa_\epsilon(\mu_{_N}[\bx],\xi_\bx,U) \big)$.   However, the computation of this  argmin in \eqref{argmina}  can be difficult in practice. 

\vspace{2mm}

We  propose a second approach which provides a more practical derivation of an approximate optimal control for the $N$-agent MDP. It will use  randomized feedback policy for the $N$-agent model, defined as a measurable function from $\Xc^N\times [0,1]^N$ into $A^N$. 
%Let us introduce the set of $K$-Lipschitz randomized feedback policies for CMKV-MDP: 
%\beqs
%\Pi_K &=& \Big\{ \mfa \in L^0(\Pc(\Xc)\times\Xc\times [0,1];A):  \\
%& & \quad \quad  \E[d_A(\mfa(\mu, x, U), \mfa(\mu, x', U))] \;  \leq \;   K  d(x,x'), \quad \forall x,x'  \in \Xc, \; \mu \in \Pc(\Xc) \Big\}.  
%\enqs
%Notice that when the state space  $\Xc$ is finite, any randomized feedback policy $\mfa$ lies in $\Pi_K$ for $K$ $=$ $1$. \red{Au fait pourquoi $K$ $=$ $1$?} 

\begin{Theorem} \label{theo2-controlNMDP}
%Fix $K$ $\in$ $\R_+$, and 
Let $\mfa_\epsilon$ be an $\epsilon$-optimal randomized feedback policy for CMKV-MDP, assumed to 
%lie in $\Pi_K$.
satisfy the regularity condition
\begin{align}\label{cond-a}
\E[d_A(\mfa_\epsilon(\mu, x, U), \mfa_\epsilon(\mu, x', U))] & \leq \;   K  d(x,x'), \quad \forall x,x'  \in \Xc, \; \mu \in \Pc(\Xc),
\end{align}
(here $U$ $\sim$ $\Uc([0,1])$) for some positive constant $K$. 
%\red{independent of $\epsilon$}.  
Consider  the randomized feedback policy in the $N$-agent model defined by 
\beqs
\bpi_r^{\mfa_\epsilon,N}(\bx,\bu)  & := & \big(\mfa_\epsilon(\mu_{_N}[\bx],x^i,u^i)\big)_{i\in\llbracket 1,N\rrbracket},  
\enqs
for  $\bx = (x^i)_{i\in\llbracket 1,N\rrbracket} \in \Xc^N, \; \bu = (u^i)_{i\in \llbracket 1,N\rrbracket} \in [0,1]^N$. Then, the randomized feedback control $\balpha^{r,\eps,N}$ $\in$ $\bAc$ defined as 
\beqs
\balpha_t^{r,\epsilon,N} &=& \bpi_r^{\mfa_\epsilon,N}(\boX_t,\bU_t), \quad t \in \N, 
\enqs
where $\{\bU_t$ $=$ $(U_t^i)_{i\in \llbracket 1,N\rrbracket}, t\in\N\}$ is a family of mutually i.i.d. uniform random variables on $[0,1]$, independent of $\Gc$, $\beps$ $=$ $\big((\eps^i_t)_{i\in \llbracket 1, N\rrbracket},\eps_t^0\big)_{t\in\N^\star}$,  
is an  $O(\epsilon + M_N^\gamma)$-optimal control for $V_N(\bx_0)$, namely:
\beqs
V_N(\bx_0) -  C(1+K)(\eps + M_N^\gamma)  & \leq & V_N^{\balpha^{r,\eps,N}}(\bx_0).   
\enqs
%namely: there exists some  positive constant $C$ (depending on $\Delta_\Xc$, $\beta$, $K_F$, $K_f$) s.t. 
%\beqs 
% V_N(\bx_0)  -  \frac{\eps}{1-\beta} - C M_N^\gamma & \leq & J_N(\bx_0,\balpha^{\eps,N}).
 %\enqs
\end{Theorem}

\vspace{1mm}

Theorem  \ref{theo2-controlNMDP} provides a simple and natural procedure to get an approximate policy for the $N$-agent MDP:  it corresponds to using an $\epsilon$-optimal  randomized feedback policy $\mfa_\epsilon$ of the CMKV-MDP,  but instead of inputting the theoretical state distribution of the McKean-Vlasov MDP in its mean-field argument, we  input the empirical state distribution of the $N$-agent MDP, and instead of inputting the McKean-Vlasov state in its state argument, we input the $N$-agent individual states, and moreover, we use a randomization by  tossing  a coin at any time and for any agent. 
Notice that the validity of this procedure requires the Lipschitz condition \eqref{cond-a}, which always holds true when the state space  $\Xc$ is finite.  Indeed, in this case, the metric on $\Xc$ is the discrete distance $d(x,x')$ $=$ $1_{x \neq x'}$, and \eqref{cond-a} is clearly satisfied with $K$ $=$ $\Delta_A$.

\section{Proof of main results} \label{sec-chaos-N}

This section is devoted to the proofs of Theorems  \ref{theo-chaosvalue}, \ref{theo1-controlNMDP}, and \ref{theo2-controlNMDP} about  rate of convergence in the  propagation of chaos between the $N$-agent MDP and the limiting conditional  McKean-Vlasov MDP.
Our approach relies on  the Bellman operators of each MDP. By proving their proximity (in a sense to be precised), we will be able to prove on the one hand the proximity of their unique fixed points, hence the convergence of the value functions, 
and on the other hand that almost optimal randomized feedback policies are directly related to the Bellman operators via the verification result, which will give the convergence of the approximate controls.

\vspace{1mm}

\subsection{Comparing the Bellman operators}

We first introduce the following useful measurable optimal permutation for the coupling of empirical measures.

\begin{Definition}[Measurable optimal permutation] \label{defpermut} 
Let $(\Yc,d)$ be a metric space.  There exists a measurable map $\sigma: (\by,\by')\in (\Yc^N)^2\rightarrow \sigma^{\by,\by'}\in\mathfrak{S}_N$ (where $\mathfrak{S}_N$ denotes the set of permutations on $\llbracket 1, N\rrbracket$) such that for all $(\by,\by')\in (\Yc^N)^2$, we have
\begin{align} \label{couplingpermut} 
\Wc_d\big(\mu_{_N}[\by],\mu_{_N}[\by'] \big) &= \;  \bd_N(\by,\by'_{\sigma^{\by,\by'}}),
%\Wc_d\big(\frac{1}{N}\sum_{i=1}^N \delta_{x^i},\frac{1}{N}\sum_{i=1}^N \delta_{x'^i} \big) &=&. \frac{1}{N}\sum_{i=1}^N d(x^i,x'^{\sigma^{\bx,\bx'}_i}).
\end{align} 
where we set $\by'_{\sigma^{\by,\by'}}$ $=$ $(y'^{\sigma^{\by,\by'}_i})_{i\in\llbracket 1,N\rrbracket}$ for $\by'$ $=$ $(y'^i)_{i\in\llbracket 1,N\rrbracket}$. 
\end{Definition}
\noindent{\bf Proof.}
It is a well known result (see \cite{thorpe2018introduction}) that, given $(\by,\by')\in (\Yc^N)^2$, there exists a permutation $\sigma^{\by,\by'}\in \mathfrak{S}_N$ realizing an optimal coupling between  $\mu_{N}[\by]$, $\mu_{_N}[\by']$ $\in$ 
$\Pc(\Yc)$, i.e., s.t. \eqref{couplingpermut} holds. 
Let us check  that this optimal permutation can be represented as a measurable function of $(\by,\by')\in (\Yc^N)^2$. Let $n\in \llbracket 1, N!\rrbracket\mapsto \sigma^n \in \mathfrak{S}_N$ be some  bijection. Notice that the function 
\beqs 
\by,\by'\in \Yc^N \mapsto \left(\bd_N(\by,\by'_{\sigma^n}) \right)_{n\in \llbracket 1, N!\rrbracket}\in \R^{N!}
\enqs 
is continuous, hence  measurable. Furthermore, it is clear that the function
\beqs 
\bz\in \R^{N!}\mapsto  \min \big[ \argmin_{n\in N!}z^n) \big] 
\enqs 
is measurable. Denoting by 
\beqs 
n_{min}(\by,\by') &:= &  \min \big[ \argmin_{n\in N!} \bd_N(\by,\by'_{\sigma^n}) \big], 
\enqs 
it follows that  the function $\by,\by'\in \Xc^N\mapsto \sigma^{\by,\by'}$ $=$ $\sigma^{n_{min}(\by,\by')}$ is a measurable representation of the optimal permutation.
\ep

\vspace{3mm}

We now study the ``proximity'' between the Bellman  operator of the CMKV-MDP given in \eqref{defTa}, and the Bellman operator of  the $N$-agent problem,  viewed  as a MDP  with state space $\Xc^N$, action space $A^N$, noise sequence $\beps=(\beps_t)_{t\in\N^\star}$ with  $\beps_t$ $:=$ $((\eps^i_t)_{i\in\llbracket 1, N\rrbracket}, \eps^0_t)$  valued in $E^N\times E^0$,  state transition function
\beqs 
\bF(\bx, \ba, \be) &:= & \Big(F(x^i, a^i, \mu_{_N}[\bx,\ba], e^i, e^0)\Big)_{i\in \llbracket 1,N\rrbracket},  \quad  \be=((e^i)_{i\in \llbracket 1,N\rrbracket}, e^0)\in E^N\times E^0,
\enqs 
and reward function 
\beqs
\bff (\bx,\ba) &= &\frac{1}{N}\sum_{i=1}^N f\big(x^i, a^i, \mu_{_N}[\bx,\ba] \big), \quad \bx = (x^i)_{i\in \llbracket 1,N\rrbracket}, \; \ba = (a^i)_{i\in \llbracket 1,N\rrbracket}. 
\enqs 
Denoting by  $L^\infty_m(\Xc^N)$ the subset of measurable functions in $L^\infty(\Xc^N)$ (the set of bounded real-valued functions on $\Xc^N$),  the Bellman ``operator'' $\Tc_N:L^\infty_m(\Xc^N)\rightarrow L^\infty(\Xc^N)$ of the $N$-agent MDP is defined 
for any $W$ $\in$ $L^\infty_m(\Xc^N)$ by: 
\begin{align} 
\Tc_N W(\bx) & := \; \sup_{\ba \in A^N} \T^{\ba}_N W (\bx),  \quad \bx \in \Xc^N, 
\end{align}
where
\begin{align}  
\T^{\ba}_N W (\bx) & := \; \bff(\bx,\ba) + \beta \E \big[ W\big(  \bF(\bx,\ba, \beps_1)\big) \big],  \quad \bx \in \Xc^N,  \; \ba \in A^N. 
\end{align}
The characterization of the value function $V_N$ and optimal controls for the $N$-agent MDP via the Bellman operator $\Tc_N$ is stated  in Appendix  \ref{sec-Bell-N}.

We aim to quantify  how ``close'' $\T^\ba_N$ and $\T^\mra$ are when $\ba$ and $\mra$ are close in a sense to be precised. 
Notice that the $N$-agent operator $\T^\ba_N$ is defined on $L^\infty_m(\Xc^N)$ while the McKean-Vlasov operator $\T^\mra$ is defined on $L^\infty_m(\Pc(\Xc))$. There is however a natural way to compare them by means of an ``unlifting'' procedure. To any function $W\in L^\infty_m(\Pc(\Xc))$, we associate the unlifted function $\widecheck{W}\in L^\infty_m(\Xc^N)$ defined by
$$
\widecheck{W}(\bx):=W(\mu_{_N}[\bx]), \quad \forall \bx \in \Xc^N.
$$
%In other words, by associating to any $\bx\in \Xc^N$ its empirical distribution $\frac{1}{N}\sum_{n=1}^N \delta_{x^n}\in\Pc(\Xc)$, we can see any function $W\in L^\infty_m(\Pc(\Xc))$ as a function $\widecheck{W}\in L^\infty_m(\Xc^N)$. 
%This trick will allow us to compare $\T^\ba_N$ and $\T^a$.
%The first step is to compare the Bellman operators. Let $a(\xi,U)$ be a randomized feedback policy, and let $\ba$ be a random action valued in $A^N$. Let us denote
%\beqs 
%\epsilon:= \E[\Wc(\P_{\xi,a(\xi,U)}, \frac{1}{N}\sum_{n=1}^N \delta_{x^n,\ba^n})]
%\enqs 
%Our goal is to compare $\Tc^{a}V^\star$ to $\Tc_N^{\ba}V^\star$.
We recall from \cite{motte2019meanfield} that the value function $V$ of the CMKV-MDP is $\gamma$-H\"older:
\begin{align}  \label{Vholder} 
| V(\mu) - V(\mu') | & \leq \; K_\star  \big( \Wc_d(\mu,\mu')  \big)^\gamma, \quad \forall \mu,\mu' \in \Pc(\Xc),  
\end{align} 
for some constant $K_\star$ depending on $K_F$, $\beta$ and $\Delta_{\Xc}$.

\begin{Lemma}\label{lem-bellman-comp}
There exists some positive constant $C$ such that for all  $\mra\in L^0(\Xc\times [0,1];A)$, $\ba \in A^N$, $\bx\in\Xc^N$ and $(\xi_\bx,U)$ $\sim$ $\mu_{_N}[\bx]\otimes \Uc([0,1])$, 
\beqs 
\vert\widecheck{\T^{\mra}V}(\bx) - \T_N^{\ba}\widecheck{V}(\bx)\vert&\leq& C \Big[ 
%\Wc_{\bomrd}(\Lc(\xi,\mra(\xi,U)), \mu_{_N}[\bx,\ba])  + 
\big(\Wc_{\bomrd}(\Lc(\xi_\bx,\mra(\xi_\bx,U)), \mu_{_N}[\bx,\ba])  \big)^\gamma + M_N^\gamma \Big]. 
\enqs 
\end{Lemma}
{\bf Proof.} For any $\mra\in L^0(\Xc\times [0,1];A)$, and $\ba$ $=$ $(a^i)_{i\in\llbracket 1,N\rrbracket}$ $\in$ $A^N$, we have
\begin{align} 
& \widecheck{\T^{\mra}V}(\bx) - \T_N^{\ba}\widecheck{V}(\bx) \nonumber \\
%&=&\E\Big[ \E \big[ f(\xi,\mra(\xi,U),\Lc(\xi,\mra(\xi,U))) + \beta V\big(  \P^0_{F(\xi,\mra(\xi,U),\Lc(\xi,\mra(\xi,U)),\eps_1,\eps_1^0)} \big) \big]. \\
%& & \quad - \;  \Big(\frac{1}{N} \sum_{i=1}^N f(x^i,a^i,\mu_{_N}[\bx,\ba] ) + \beta V\big(  \frac{1}{N}\sum_{j=1}^N \delta_{F(x^j,a^j, \mu_{_N}[\bx,\ba], \eps_1,\eps_1^0)}   \big)\Big) \Big]\\
&= \; \E\Big[ f(\xi_\bx,a(\xi_\bx,U),\Lc(\xi_\bx,\mra(\xi_\bx,U)) )  -  \frac{1}{N} \sum_{i=1}^N f(x^i,a^i,\mu_{_N}[\bx,\ba]) \Big]  \nonumber \\
&   \quad + \; \beta \E\Big[  V\big(  \P^0_{F(\xi_\bx,\mra(\xi_\bx,U),\Lc(\xi_\bx,\mra(\xi_\bx,U)),\eps_1,\eps_1^0)}   \big)  - V\big(  \frac{1}{N}\sum_{i=1}^N \delta_{F(x^i,a^i, \mu_{_N}[\bx,\ba],\eps_1^i,\eps_1^0)}   \big)  \Big]. \label{decomp} 
\end{align} 
We write 
\beqs 
\E \Big[ f(\xi_\bx,a(\xi_\bx,U),\Lc(\xi_\bx,\mra(\xi_\bx,U)))\Big] -\frac{1}{N} \sum_{i=1}^N f(x^i,\ba^i,\mu_{_N}[\bx,\ba]) &=& \hat{f}(\Lc(\xi_\bx,\mra(\xi_\bx,U))) -  \hat{f}(\mu_{_N}[\bx,\ba]),
\enqs 
where $\hat{f}(\mu)$ $=$ $\int f(x',a',\mu) \mu(\d x',\d a')$ for all $\mu\in \Pc(\Xc\times A)$. Notice that for $\mu,\mu'\in \Pc(\Xc\times A)$, we have
\beqs 
\hat{f}(\mu)-\hat{f}(\mu')&=&\int f(x',a',\mu)(\mu-\mu')(\d x',\d a') +\int (f(x',a',\mu)-f(x',a',\mu'))\mu'(\d x',\d a')\\
&\leq&K_f\Wc_{\bomrd}(\mu,\mu')+ K_f\Wc_{\bomrd}(\mu,\mu') \;= \; 2K_f\Wc_{\bomrd}(\mu,\mu'),
\enqs 
from the Kantorovich-Rubinstein dual representation \eqref{dualW} and $({\bf Hf_{lip}})$. It follows that 
\begin{align} \label{eq-flip}
& \Big|  \E \Big[ f\big(\xi_\bx,\mra(\xi_\bx,U),\Lc(\xi_\bx,\mra(\xi_\bx,U))\big)\Big] -\frac{1}{N} \sum_{i=1}^N f(x^i,a^i,\mu_{_N}[\bx,\ba]) \Big| \\
\; \leq & \;\;  2K_f\Wc_\bomrd(\Lc(\xi_\bx,\mra(\xi_\bx,U)), \mu_{_N}[\bx,\ba]).
\end{align}
Let us next  focus on the second term in \eqref{decomp}.  
%\beqs 
%\E\Big[ \E \Big[V\big(  \P^0_{F(\xi,a(\xi,U),\P_{\xi,a(\xi,U)},\eps^i_1,\eps_1^0}   \big) \Big] - V\big(  \frac{1}{N}\sum_{i=1}^N \delta_{F(x^i,a^i, \frac{1}{N}\sum_{j=1}^N \delta_{x^j,a^j},\eps^i_1,\eps_1^0}   \big)\Big) \Big]
%\enqs 
As $V$ is $\gamma$-Hölder with constant factor $K_\star$, we have
\begin{align} 
&  \Big|   \E\Big[  V\big(  \P^0_{F(\xi_\bx,\mra(\xi_\bx,U),\Lc(\xi_\bx,\mra(\xi_\bx,U)),\eps_1,\eps_1^0)}   \big)  - V\big(  \frac{1}{N}\sum_{i=1}^N \delta_{F(x^i,a^i, \mu_{_N}[\bx,\ba],\eps_1^i,\eps_1^0)}   \big)  \Big]  \Big| \nonumber \\
%&\leq&K_\star\E \Big[\Wc_d( \P^0_{F(\xi,a(\xi,U),\Lc(\xi,a(\xi,U)),\eps_1,\eps_1^0)} ,\frac{1}{N}\sum_{i=1}^N \delta_{F(x^i,a^i, \frac{1}{N}\sum_{j=1}^N \delta_{x^j,a^j},\eps^i_1,\eps_1^0)}   \big)^\gamma \Big]\\
\leq & \; K_\star\E \Big[\Wc_d( \P^0_{F(\xi_\bx,\mra(\xi_\bx,U),\Lc(\xi_\bx,\mra(\xi_\bx,U)),\eps^i_1,\eps_1^0)} ,\frac{1}{N}\sum_{i=1}^N \delta_{F(x^i,a^i,  \mu_{_N}[\bx,\ba],\eps^i_1,\eps_1^0)}   \big) \Big]^\gamma, \label{interV0} 
\end{align}  
by  Jensen's inequality. Let $(\xi^i,U_0^i)_{i\in\llbracket 1, N\rrbracket}$ be $N$ i.i.d. random variables, independent of $\beps_1$, such that $(\xi^i,U_0^i)$ $\sim$ $\mu_{_N}[\bx]\otimes\Uc([0,1])$, $i\in\llbracket 1, N\rrbracket$. For any i.i.d. random variables  $(\tilde{\eps}^i_1)_{i\in \llbracket 1, N\rrbracket}$ such that
\begin{align}\label{eq-eps-tilde}
    ((\xi^i, U^i_0,\tilde{\eps}^i_1)_{i\in \llbracket 1, N\rrbracket},\eps^0_1)\overset{d}{=} ((\xi^i, U^i_0,\eps^i_1)_{i\in \llbracket 1, N\rrbracket},\eps^0_1),
\end{align}
we have 
\begin{align}  
&  \E \Big[\Wc_d\big( \P^0_{F(\xi_\bx,\mra(\xi_\bx,U),\Lc(\xi_\bx,\mra(\xi_\bx,U)),\eps^i_1,\eps_1^0)} ,\frac{1}{N}\sum_{i=1}^N \delta_{F(x^i,a^i,  \mu_{_N}[\bx,\ba],\eps^i_1,\eps_1^0)}   \big) \Big] \nonumber \\
&\leq \; \E\Big[\Wc_d\big(  \P^0_{F(\xi_\bx,\mra(\xi_\bx,U),\Lc(\xi_\bx,\mra(\xi_\bx,U)),\eps^i_1,\eps_1^0)}  ,   \frac{1}{N}\sum_{i=1}^N \delta_{F(\xi^i,\mra(\xi^i,U^i_0),\Lc(\xi_\bx,\mra(\xi_\bx,U)),\tilde{\eps}^i_1,\eps_1^0)}  \big) \Big] \nonumber \\
&  \; + \; \E\Big[\Wc_d\big( \frac{1}{N}\sum_{i=1}^N \delta_{F(\xi^i,\mra(\xi^i,U_0^i),\Lc(\xi_\bx,\mra(\xi_\bx,U)),\tilde{\eps}^i_1,\eps_1^0)} ,   \frac{1}{N}\sum_{i=1}^N \delta_{F(x^i,a^i,  \mu_{_N}[\bx,\ba],\eps^i_1,\eps_1^0}  \big) \Big] \nonumber \\
&\leq \; M_N+\E\Big[\Wc_d\big( \frac{1}{N}\sum_{i=1}^N \delta_{F(\xi^i,\mra(\xi^i,U_0^i),\Lc(\xi_\bx,\mra(\xi_\bx,U)),\tilde{\eps}^i_1,\eps_1^0)} ,   \frac{1}{N}\sum_{i=1}^N \delta_{F(x^i,a^i,  \mu_{_N}[\bx,\ba],\eps^i_1,\eps_1^0}  \big) \Big],  \label{interW1} 
\end{align} 
by definition of $M_N$ in \eqref{defMN}. 
%Let us now focus on 
%\beqs 
%\E\Big[\Wc\left( \frac{1}{N}\sum_{i=1}^N \delta_{F(\xi^i,a(\xi^i,U^i_1),\Lc(\xi,a(\xi,U)),\tilde{\eps}^i_1,\eps_1^0)} ,   \frac{1}{N}\sum_{i=1}^N \delta_{F(x^i,a^i, \frac{1}{N}\sum_{j=1}^N \delta_{x^j,a^j},\eps^i_1,\eps_1^0}  \right) \Big].
%\enqs 
%The reason why we considered an arbitrary random family $(\tilde{\eps}^i_1)_{i\in \llbracket 1, N\rrbracket}$ with same distribution as $(\eps^i_1)_{i\in \llbracket 1, N\rrbracket}$ instead of just considering $(\eps^i_1)_{n\in \llbracket 1, N\rrbracket}$ was to allow us to {\it couple} variables nicely before using the formula 
%$\Wc(\frac{1}{N}\sum_{i=1}^N \delta_{y^i},\frac{1}{N}\sum_{i=1}^N \delta_{z^i})\leq \frac{1}{N}\sum_{i=1}^N d(y^i,z^i)$ in order to obtain a good estimation. 
Let us now  consider the random permutation $\sigma^{(\xi^i,\mra(\xi^i,U^i_0))_{i\in\llbracket 1,N\rrbracket}, (x^i,a^i)_{i\in\llbracket 1, N\rrbracket}}$ defined in Definition \ref{defpermut} that we shall, to simplify notations, simply denote by $\sigma$. 
Notice that as $(\xi^i,\mra(\xi^i,U^i_0))_{i\in\llbracket 1, N\rrbracket}\independent (\eps^i_1)_{i\in\llbracket 1, N\rrbracket}$, we clearly see that 
$(\tilde{\eps}^i_1)_{i\in\llbracket 1, N\rrbracket}:=(\eps^{(\sigma^{-1})_i}_1)_{i\in\llbracket 1, N\rrbracket}$ satisfies the required condition \eqref{eq-eps-tilde}. Therefore the above relation applies to $(\tilde{\eps}^i_1)_{i\in\llbracket 1, N\rrbracket}= (\eps^{(\sigma^{-1})_i}_1)_{i\in\llbracket 1, N\rrbracket}$. 
For such $(\tilde{\eps}^i_1)_{i\in\llbracket 1, N\rrbracket}$, we get 
\beqs 
&&\E\Big[\Wc_d\big( \frac{1}{N}\sum_{i=1}^N \delta_{F(\xi^i,\mra(\xi^i,U^i_0),\Lc(\xi_\bx,\mra(\xi_\bx,U)),\eps^{(\sigma^{-1})_i}_1,\eps_1^0)} ,   \frac{1}{N}\sum_{i=1}^N \delta_{F(x^i,a^i,\mu_{_N}[\bx,\ba],\eps^i_1,\eps_1^0)}  \big) \Big]\\
&=&\E\Big[\Wc_d\big( \frac{1}{N}\sum_{i=1}^N \delta_{F(\xi^{\sigma_i},\mra(\xi^{\sigma_i},U^{\sigma_i}_0),\Lc(\xi_\bx,\mra(\xi_\bx,U)),\eps^{i}_1,\eps_1^0)} ,   \frac{1}{N}\sum_{i=1}^N \delta_{F(x^i,a^i,\mu_{_N}[\bx,\ba],\eps^i_1,\eps_1^0)}  \big) \Big]\\
&\leq&\frac{1}{N}\sum_{i=1}^N\E\Big[ d\big(F(\xi^{\sigma_i},\mra(\xi^{\sigma_i},U^{\sigma_i}_0),\Lc(\xi_\bx,\mra(\xi_\bx,U)),\eps^{i}_1,\eps_1^0) ,   F(x^i,a^i,\mu_{_N}[\bx,\ba],\eps^i_1,\eps_1^0)\big)  \Big] \\
%\enqs 
%By conditioning w.r.t. $((\xi_i,U_i)_{i\in\llbracket 1, N\rrbracket},\eps^0_1)$ and using the regularity in expectation of $F$ given by  $({\bf HF_{lip}})$, we obtain
%\beqs 
%&&\frac{1}{N}\sum_{i=1}^N\E\Big[ d(F(\xi^{\sigma_i},a(\xi^{\sigma_i},U^{\sigma_i}_1),\P_{\xi,a(\xi,U)},\eps^{i}_1,\eps_1^0) ,   F(x^i,a^i,\frac{1}{N}\sum_{j=1}^N \delta_{x^j,a^j},\eps^i_1,\eps_1^0))  \Big]\\
&\leq& K_F\frac{1}{N}\sum_{i=1}^N\E\Big[ \bomrd\big((\xi^{\sigma_i},\mra(\xi^{\sigma_i},U^{\sigma_i}_0)),(x^i,a^i) \big)  + \Wc_{\bomrd}\big(\Lc(\xi_\bx,\mra(\xi_\bx,U)),\mu_{_N}[\bx,\ba]\big) \Big]\\
&=& K_F\E[\Wc_{\bomrd}(\frac{1}{N}\sum_{i=1}^N \delta_{(\xi^i,\mra(\xi^{i},U^{i}_0))},\mu_{_N}[\bx,\ba])+ \Wc_{\bomrd}\big(\Lc(\xi_\bx,\mra(\xi_\bx,U)),\mu_{_N}[\bx,\ba] \big)  \Big]\\
&\leq& K_F \Big(M_N+2\E\big[\Wc_{\bomrd}(\Lc(\xi_\bx,\mra(\xi_\bx,U)), \mu_{_N}[\bx,\ba] \big) \big] \Big),
\enqs 
where the first inequality comes from \eqref{inegW}, the second one is derived by conditioning w.r.t. $((\xi^i,U_0^i)_{i\in\llbracket 1, N\rrbracket},\eps^0_1)$ and using the regularity in expectation of $F$ in $({\bf HF_{lip}})$,  the 
last  equality holds true by definition of the permutation $\sigma$ realizing the optimal coupling \eqref{couplingpermut},  and the last inequality from the definition of $M_N$. Recalling \eqref{interW1}, we then have  
\beqs 
&& \E \Big[\Wc_d\big( \P^0_{F(\xi_\bx,\mra(\xi_\bx,U),\Lc(\xi_\bx,\mra(\xi_\bx,U)),\eps^i_1,\eps_1^0)} ,\frac{1}{N}\sum_{i=1}^N \delta_{F(x^i,a^i,  \mu_{_N}[\bx,\ba],\eps^i_1,\eps_1^0)}   \big) \Big]  \\
%&\leq& M_N+K_F(M_N+2\E[\Wc(\P_{\xi,a(\xi,U)}, \frac{1}{N}\sum_{i=1}^N \delta_{x^i,a^i})])\\
&\leq& (1+K_F)M_N+2K_F\E\big[\Wc_{\bomrd}(\Lc(\xi_\bx,\mra(\xi_\bx,U)), \mu_{_N}[\bx,\ba] \big) \big] 
\enqs 
which implies by \eqref{interV0} 
\beqs 
&&  \E\Big[  V\big(  \P^0_{F(\xi_\bx,\mra(\xi_\bx,U),\Lc(\xi,\mra(\xi,U)),\eps_1,\eps_1^0)}   \big)  - V\big(  \frac{1}{N}\sum_{i=1}^N \delta_{F(x^i,a^i, \mu_{_N}[\bx,\ba],\eps_1^i,\eps_1^0)}   \big)  \Big]   \\
&\leq&K_\star \Big( (1+K_F)M_N+2K_F\E\big[\Wc_{\bomrd}(\Lc(\xi_\bx,\mra(\xi_\bx,U)), \mu_{_N}[\bx,\ba] \big) \big]  \Big)^\gamma. 
\enqs 
Together with \eqref{eq-flip}, and plugging into \eqref{decomp}, we obtain finally 
\beqs 
&& \Big| \widecheck{\T^{a}V}(\bx) - \T_N^{\ba}\widecheck{V}(\bx) \Big| \\
&\leq&  2K_f  \E \Big[  \Wc_\bomrd(\Lc(\xi_\bx,\mra(\xi_\bx,U)), \mu_{_N}[\bx,\ba]) \Big]  + K_\star \Big( (1+K_F)M_N+2K_F\E\big[\Wc_{\bomrd}(\Lc(\xi_\bx,\mra(\xi_\bx,U)), \mu_{_N}[\bx,\ba] \big) \big]  \Big)^\gamma \\
&\leq&  C  \Big\{ \Wc_\bomrd(\Lc(\xi_\bx,\mra(\xi_\bx,U)), \mu_{_N}[\bx,\ba])  +   \Big( \Wc_\bomrd(\Lc(\xi_\bx,\mra(\xi_\bx,U)), \mu_{_N}[\bx,\ba]) \Big)^\gamma + M_N^\gamma \Big\}
%\Oc(\E[\Wc(\P_{\xi,a(\xi,U)}, \frac{1}{N}\sum_{i=1}^N \delta_{x^i,a^i})]^\gamma+M_N^\gamma).
\enqs
(recall that $\gamma$ $\leq$ $1$), for some constant $C$ depending only on $K_\star$, $K_f$, $K_F$,  where we also use the fact that $\Wc_\bomrd(\Lc(\xi_\bx,\mra(\xi_\bx,U)), \mu_{_N}[\bx,\ba])$ is bounded by a constant depending on the diameter of the compact set $\Xc\times A$.  
This ends the proof.
\ep

\vspace{3mm}

%Lemma \ref{lem-bellman-comp} provides a general comparison between the operators of the $N$-agent and the McKean-Vlasov MDPs. Let us now exploit this result to derive more specific and interesting links between both problems. In Section \ref{sec-chaos-value}, we establish the proximity between their value functions, and in Section \ref{sec-%chaos-policy}, we provide ways to derive a good policy for the $N$-agent MDP from a good policy for the McKean-Vlasov MDP.

\subsection{Proof of Theorem \ref{theo-chaosvalue} }\label{sec-chaos-value}

Lemma \ref{lem-bellman-comp} means  that given $\mra$ $\in$ $L^0(\Xc\times [0,1];A)$, $\ba\in A^N$, and for $\bx$ $\in$ $\Xc^N$, the Wasserstein distance between the distribution law of $(\xi_\bx, \mra(\xi_\bx, U))$  (where $(\xi_\bx,U)$ $\sim$ $\mu_{_N}[\bx]\otimes\Uc([0,1])$),  and the empirical measure 
$\mu_{_N}[\bx,\ba]$  is small (and $N$ large),  then $\T^aV$ $\simeq$  $\T^{\ba}_N \widecheck{V}$. It is thus natural to look for suitable choices of $\mra$ $\in$ $L^0(\Xc\times [0,1];A)$, $\ba\in A^N$ so that  the above Wasserstein distance is as small as possible.  This is quantified  in the  following result.

\begin{Lemma}\label{lem-action-couple}
Fix $\bx$ $\in$ $\Xc^N$. Then, for any $\mra$ $\in$ $L^0(\Xc\times [0,1];A)$,  there exists $\ba^{\mra}\in A^N$ 
%(\red{depend de $\bx$?}\bl{Oui dépend de $\bx$. On peut soit fixer $\bx$ pour toute la section, soit écrire $\ba^{\bx, \mra}$ et $\mra^{\bx,\ba}$}) 
such that 
\beqs 
\Wc_\bomrd\big(\Lc(\xi_\bx, \mra(\xi_\bx, U)), \mu_{_N}[\bx,\ba^{\mra}]\big) & \leq & 2 M_N,
\enqs 
where $(\xi_\bx,U)$ $\sim$ $\mu_{_N}[\bx]\otimes\Uc([0,1])$. 
Conversely, for any $\ba\in A^N$, there exists $\mra^{\ba}\in L^0(\Xc\times [0,1];A)$ such that
\beqs 
\Wc_\bomrd\big(\Lc(\xi_\bx, \mra^\ba(\xi_\bx, U)), \mu_{_N}[\bx,\ba^{}]\big)  &=& 0.
\enqs 
\end{Lemma}
{\bf Proof.}
Fix $\mra\in L^0(\Xc\times [0,1];A)$. Let us consider $\bxi=(\xi^i)_{i\in \llbracket 1, N\rrbracket}$ i.i.d. with common distribution $\mu_{_N}[\bx],$ independent from $\bU_0$ $=$ $(U_0^i)_{i\in \llbracket 1, N\rrbracket}$ i.i.d. $\sim$ $\Uc([0,1])$.  We have
\beqs 
&& \E\big[\Wc_\bomrd\big(\Lc(\xi_\bx, \mra(\xi_\bx, U)),  \frac{1}{N}\sum_{i=1}^N \delta_{x^i,\mra(\xi^{\sigma^{\bxi,\bx}_i}, U^i_0)} \big)\big] \\
&\leq& \E\big[\Wc\big(\Lc(\xi_\bx,\mra(\xi_\bx, U)),\frac{1}{N}\sum_{i=1}^N \delta_{\xi^{\sigma^{\bxi,\bx}_i}, \mra(\xi^{\sigma^{\bxi,\bx}_i}, U^i_0)}\big) 
+ \Wc_{\bomrd}(\frac{1}{N}\sum_{i=1}^N \delta_{\xi^{\sigma^{\bxi,\bx}_i}, \mra(\xi^{\sigma^{\bxi,\bx}_i}, U^i_0)}, \frac{1}{N}\sum_{i=1}^N \delta_{x^i,\mra(\xi^{\sigma^{\bxi,\bx}_i}, U^i_0)}\big) \big] \\
%&\leq& M_N+\E[\Wc(\frac{1}{N}\sum_{i=1}^N \delta_{\xi^{\sigma^{\bxi,\bx}_i}, a(\xi^{\sigma^{\bxi,\bx}_i}, U^i_1)}, \frac{1}{N}\sum_{i=1}^N \delta_{x^i,a(\xi^{\sigma^{\bxi,\bx}_i}, U^i_1)})]  \\
&\leq&M_N+\E[\frac{1}{N}\sum_{i=1}^Nd(\xi^{\sigma^{\bxi,\bx}_i},x^i)] \; \leq \;  2M_N,
\enqs 
where we used the definition of $M_N$ and  \eqref{inegW} in the second inequality, and definition of $\sigma^{\bxi,\bx}$ in the last inequality. It follows that
\beqs 
\P\big[ \Wc_\bomrd\big(\Lc(\xi_\bx, \mra(\xi_\bx, U)),  \frac{1}{N}\sum_{i=1}^N \delta_{x^i,\mra(\xi^{\sigma^{\bxi,\bx}_i}, U^i_0)} \big) \leq 2M_N \big] & > & 0,
\enqs 
which implies that there exists a vector $\ba\in A^N$ 
%(\red{pourquoi, est ce clair?}\bl{Comme l'événement a une proba strictement positive, il existe un $\omega\in\Omega$ tel que
%\beqs 
%\Wc_\bomrd\big(\Lc(\xi_\bx, \mra(\xi_\bx, U)),  \frac{1}{N}\sum_{i=1}^N \delta_{x^i,\mra(\xi^{\sigma^{\bxi(\omega),\bx}_i}(\omega), U^i_0(\omega))} \big) \leq 2M_N
%\enqs
%On prend alors $\ba=(\mra(\xi^{\sigma^{\bxi(\omega),\bx}_i}(\omega), U^i_1(\omega)))_{i\in\llbracket 1, N\rrbracket}$.}  
such that
\beqs 
\Wc_\bomrd\big(\Lc(\xi_\bx, \mra(\xi_\bx, U)), \mu_{_N}[\bx,\ba]\big) &\leq&  2M_N.
\enqs 
On the other hand, given such an $\ba\in A^N$, there clearly exists $\mra^\ba\in L^0(\Xc\times [0,1];A)$ such that $\Lc(\xi_\bx,\mra^\ba(\xi_\bx,U))$ $=$ $\mu_{_N}[\bx,\ba]$: indeed, by considering $(\tilde{\xi}, \tilde{\alpha})$ $\sim$ $\mu_{_N}[\bx,\ba]$, it suffices to choose $\mra^\ba$ as a kernel for simulating the conditional distribution of $\tilde{\alpha}$ knowing 
$\tilde{\xi}$. We then have
\beqs 
\Wc\big(\Lc(\xi_\bx, \mra^\ba(\xi_\bx, U)), \mu_{_N}[\bx,\ba] \big) &=& 0.
\enqs 
\ep

\vspace{3mm}

By combining the general comparison of Bellman operators in Lemma \ref{lem-bellman-comp} with the coupling result in Lemma \ref{lem-action-couple}, we can now prove the propagation of chaos of value functions.

\vspace{3mm}

\noindent  {\bf Proof of Theorem \ref{theo-chaosvalue}.} 
From the fixed point equation for $V$ with Bellman operator $\Tc$ in \eqref{defTa}, we have
\beqs 
\widecheck{V}(\bx)&=&\widecheck{\Tc V} (\bx)\\
&=&\sup_{\mra\in L^0(\Xc\times [0,1];A)}\widecheck{\T^a V} (\bx) \; \leq \;  \sup_{\mra\in L^0(\Xc\times [0,1];A)}\T^{\ba^{\mra}}_N \widecheck{V} (\bx)+ C M_N^\gamma  \\
&\leq&\Tc_N \widecheck{V} (\bx)+  C M_N^\gamma,
\enqs 
where we used Lemma \ref{lem-bellman-comp} and Lemma \ref{lem-action-couple} in the first  inequality, and the definition of $\Tc_N$ in the last one. 
Since  $V_N$ is a fixed point of $\Tc_N$ (see Proposition \ref{theomainDPP-N}),  we then have:
\beqs 
(\widecheck{V} - V_N)(\bx)&\leq& (\Tc_N \widecheck{V} - \Tc_N V_N)(\bx) +  C M_N^\gamma, 
\enqs 
and thus by definition of $\Tc_N$, 
\beqs 
(\widecheck{V}- V_N)(\bx)&\leq& \beta \sup_{\bx'\in\Xc^N}(\widecheck{V}- V_N)(\bx')+ C M_N^\gamma, 
\enqs 
which implies 
\beqs 
\sup_{\bx\in\Xc^N}(\widecheck{V}-V_N)(\bx)&\leq& C M_N^\gamma. 
\enqs 
Likewise, by Lemma \ref{lem-bellman-comp} and Lemma \ref{lem-action-couple}, we have 
\beqs 
\widecheck{V}(\bx)&=&\widecheck{\Tc V} (\bx) \; = \; \sup_{\mra\in L^0(\Xc\times [0,1];A)}\widecheck{\T^{\mra} V} (\bx) \; \geq \;  \sup_{\ba\in A^N}\widecheck{\T^{\mra^\ba} V} (x)\\
&\geq& \sup_{\ba\in A^N}\T^{\ba}_N \widecheck{V} (\bx) -  C M_N^\gamma \; =  \; \Tc_N \widecheck{V} (\bx)-  C M_N^\gamma, 
\enqs
and using the fact that  $V_N$ is a fixed point of $\Tc_N$, we obtain similarly
%\beqs 
%(V_N-\widecheck{V})(\bx)&\leq& (\Tc_N V_N-\Tc_N \widecheck{V})(\bx) +  C M_N^\gamma, 
%\enqs 
%and thus by definition of $\Tc_N$, 
%\beqs 
%(V_N-\widecheck{V})(\bx)&\leq& \beta \sup_{\bx'\in\Xc^N}(V_N-\widecheck{V})(\bx')+ C M_N^\gamma, 
%\enqs 
%which implies 
\beqs 
\sup_{\bx\in\Xc^N}(V_N-\widecheck{V})(\bx)&\leq& C M_N^\gamma,
\enqs 
%and thus $V_N(\bx)\leq \widecheck{V}(\bx)+\Oc(M_N^\gamma)$ for all $\bx\in\Xc^N$. 
%By a similar argument, one can prove that $V_N(\bx)\geq \widecheck{V}(\bx)+\Oc(M_N^\gamma)$, and thus, $\Vert \widecheck{V}-V_N\Vert= \Oc(M_N^\gamma)$, 
which concludes the proof.

\subsection{Proof of Theorem \ref{theo1-controlNMDP} }\label{sec-chaos-policy1}

%In this section, we provide two ways to derive a good policy for the $N$-agent MDP from a good policy for the McKean-Vlasov MDP. We first provide a general approach, and then another approach more practical but relying on some additional assumption.

%\subsubsection{General case: a natural $N$-agent feedback policy}

We start with a general result estimating the efficiency of a feedback policy for the $N$-agent MDP by ``comparing'' it to an  $\epsilon$-optimal  randomized feedback policy for the CMKV- MDP.

\begin{Lemma}\label{lem-verif-res}
Let $\mfa_\epsilon$ be an $\epsilon$-optimal randomized feedback policy for the CMKV-MDP, and $\ba\in A^N$. Then, there exists some positive constant $C$ (depending only  on $\Delta_{\Xc\times A}$, $\beta$, $K_F$, $K_f$) such that for all $\bx$ $\in$ $\Xc^N$, 
\beqs 
\T^{\ba}_N V_N(\bx) &\geq&  V_N(\bx)-\epsilon- C \big[  \Wc_\bomrd\big(\Lc(\xi_\bx,\mfa_\epsilon(\mu_{_N}[\bx], \xi_\bx,U)), \mu_{_N}[\bx,\ba]\big)^\gamma + M_N^\gamma \big]. 
\enqs 
\end{Lemma}
{\bf Proof.} Fix $\bx$ $\in$ $\Xc^N$, $\ba$ $\in$ $A^N$, and  define $\mra_\epsilon$ $\in$ $L^0(\Xc\times [0,1];A)$ by $\mra_\epsilon(x,u)$ $=$ $\mfa_\epsilon(\mu_{_N}[\bx],x,u)$ for $x\in\Xc$, $u\in [0,1]$. By Theorem \ref{theo-chaosvalue} and the $\beta$-contracting property of $\T^{\ba}_N$, 
we have 
% (\red{je ne comprends pas où on utilise la $\beta$-contraction?}\bl{On a
\beqs 
\vert\T^{\ba}_NV_N(\bx)-\T_N^{\ba}\widecheck{V}(\bx)\vert &\leq& \beta \Vert V_N(\bx) -\widecheck{V}(\bx) \Vert_{\Xc^N}\ \; \leq \;   C M_N^\gamma,
\enqs
and so 
%La $\beta$-contraction est utilisée dans la première inégalité.}, we have
\beqs 
\T^{\ba}_NV_N(\bx) &\geq&  \T_N^{\ba}\widecheck{V}(\bx) -  C M_N^\gamma.
\enqs 
Together with  Lemma \ref{lem-bellman-comp}, this yields 
\begin{align}\label{eq-verif-res-0}
 %\T_N^{\ba}\widecheck{V}(\bx) 
 \T^{\ba}_NV_N(\bx) 
&\geq \;   \widecheck{\T^{\mra_\epsilon}V}(\bx)- C \big[ \Wc_\bomrd\big(\Lc(\xi,\mra_\epsilon(\xi,U)), \mu_{_N}[\bx,\ba] \big)^\gamma +  M_N^\gamma \big].
\end{align}
Denote by $\alpha^\epsilon$ the randomized feedback control associated via \eqref{ranfeedbackeps} to the randomized feedback policy $\mfa_\epsilon$. Then, notice that the gain functional $V^{\alpha^\epsilon}(\xi)$ depends on $\xi$ only through its law $\mu$ $=$ $\Lc(\xi)$, and we set  
$V^{\alpha^\epsilon}(\mu)$ $=$ $V^{\alpha^\epsilon}(\xi)$ when $\xi$ $\sim$ $\mu$.   Since $V$ $\geq$ $V^{\alpha^\epsilon}$,  and by the monotonicity of $\T^{\mra_\epsilon}$, we have 
\beqs
\T^{\mra_\epsilon} V(\mu_{_N}[\bx]) & \geq & \T^{\mra_\epsilon} V^{\alpha_\epsilon} (\mu_{_N}[\bx]) \; = \; V^{\alpha^\epsilon}(\mu_{_N}[\bx])  \; \geq \; V(\mu_{_N}[\bx]) - C \epsilon
\enqs
by recalling that $V^{\alpha^\epsilon}$ is a fixed point of $\T^{\mra_\epsilon}$, and using the fact that  $\mfa_\epsilon$ is an $\epsilon$-optimal randomized feedback policy for the CMKV-MDP.  
From Theorem \ref{theo-chaosvalue}, this implies that 
\beqs 
\widecheck{\T^{\mra_\epsilon}V}(\bx) & \geq&  V_N(\bx)- C\big[ \epsilon +  M_N^\gamma \big],
\enqs 
which proved the required result when combined with \eqref{eq-verif-res-0}. 
\ep

\vspace{3mm}

Let us denote by $L^0(\Xc^N;A^N)$ the set of measurable functions from $\Xc^N$ into $A^N$. Given a  feedback policy $\bpi$ $\in$ $L^0(\Xc^N;A^N)$ for the $N$-agent problem,  the associated feedback control is the unique control $\balpha^{\bpi}$ defined by  
$\balpha^{\bpi}_t$ $=$ $\bpi(\bX_t)$,  $t\in\N$. By misuse of notation, we denote $V^{\bpi}_N=V^{\balpha^{\bpi}}_N$. Let us then introduce the operator $\Tc^{\bpi}_N$  on $L^\infty_m(\Xc^N)$, defined  by
\begin{align}  
\Tc^{\bpi}_N W(\bx) &= \;  \bff( \bx,\bpi(\bx)) + \beta \E \big[ W\big( \bF(\bx,\bpi(\bx),\beps_1)\big) \big], \quad \bx \in \Xc^N. 
\end{align}

\begin{Proposition}\label{prop-policy-1}
Let $\mfa_\epsilon$ be an $\epsilon$-optimal randomized feedback policy for the CMKV-MDP, and  consider any  feedback policy  $\bpi$ for the $N$-agent MDP.  Then, the feedback control $\balpha^{\bpi}$ is 
\beqs 
\Oc(\epsilon + \sup_{\bx\in\Xc^N}\Wc_\bomrd\big(\Lc(\xi_\bx,\mfa_\epsilon(\mu_{_N}[\bx],\xi_\bx,U)), \mu_{_N}[\bx,\bpi(\bx)] \big)^\gamma + M_N^\gamma) \text{-optimal for } V_N(\bx_0),
\enqs 
where $(\xi_\bx,U)$ $\sim$ $\mu_{_N}[\bx]\otimes\Uc([0,1])$, namely
\beqs
V_N(\bx_0)  -  C \big[ \epsilon +   \sup_{\bx\in\Xc^N}\Wc_\bomrd\big(\Lc(\xi_\bx,\mfa_\epsilon(\mu_{_N}[\bx],\xi_\bx,U)), \mu_{_N}[\bx,\bpi(\bx)] \big)^\gamma + M_N^\gamma \big] & \leq &  V_N^{\bpi}(\bx_0). 
\enqs
\end{Proposition}
{\bf Proof.} Fix $\bx\in\Xc^N$, and let  $\ba$ $=$ $\bpi(\bx)$ $\in$ $A^N$. 
%and let $\mra_\eps$ $\in$ $L^0(\Xc\times [0,1];A)$ defined by $\mra_\eps(x,u)$ $=$ $\mfa_\eps(\mu_{_N}[\bx],x,u)$ for $x\in\Xc$, $u\in [0,1]$. 
By definition, we have  $\Tc_N^{\bpi}V_N(\bx)$ $=$ $\T^{\ba}_NV_N(\bx)$.  By Lemma \ref{lem-verif-res}, we thus have
\beqs 
\Tc^{\bpi}_NV_N(\bx)&\geq& V_N(\bx)-\epsilon- C \big[ \sup_{\bx \in \Xc^N}  \Wc_\bomrd\big(\Lc(\xi_\bx,\mfa_\epsilon(\mu_{_N}[\bx],\xi_\bx,U)), \mu_{_N}[\bx,\ba] \big)^\gamma + M_N^\gamma \big],
%\\ &\geq& V_N(\bx)-\eps-\Oc(\sup_{\bx\in\Xc^N}\Wc(\P_{\xi,a(\xi,U)}, \frac{1}{N}\sum_{i=1}^N \delta_{x^i,\ba^i})^\gamma + M_N^\gamma)
\enqs 
and we conclude by the verification result in Lemma \ref{lemverif-N}. 
\ep

\vspace{3mm}

Proposition \ref{prop-policy-1} has an important implication:  it means that a feedback policy $\bpi$ for the $N$-agent MDP yields the better performance whenever it assigns for each state $\bx$ the action $\bpi(\bx)$ that achieves the minimum of  
\beqs
\ba \in A^N & \mapsto &  \Wc_\bomrd\big(\Lc(\xi_\bx,\mfa_\eps(\mu_{_N}[\bx],\xi_\bx,U)), \mu_{_N}[\bx,\ba] \big). 
\enqs
Let us check that  one can choose a measurable version   of this argmin.

\begin{Lemma} \label{remmesur} 
Let $\mfa$ $\in$ $L^0(\Pc(\Xc)\times \Xc\times [0,1];A)$. Then, there exists a measurable function $\bpi^\star$ $:$ $\Xc^N$ $\rightarrow$ $A^N$ such that
\beqs
\bpi^\star(\bx) & \in & \argmin_{\ba\in A^N}  \Wc_\bomrd\big(\Lc(\xi_\bx,\mfa(\mu_{_N}[\bx],\xi_\bx,U)), \mu_{_N}[\bx,\ba] \big), \quad \bx \in \Xc^N. 
\enqs
\end{Lemma}
{\bf Proof.} 
Notice that the function
\beqs 
h(\bx,\ba) &:=&  \Wc_\bomrd\big(\Lc(\xi_\bx,\mfa(\mu_{_N}[\bx],\xi_\bx,U)), \mu_{_N}[\bx,\ba] \big)
\enqs 
is such that for $\ba\in A^N$, $h(\cdot, \ba)$ is measurable, and for $\bx\in\Xc^N$, $h(\bx, \cdot)$ is continuous. Let us then show that  one can measurably select $\argmin_{\ba\in A^N} h(\bx, \ba)$ w.r.t. $\bx$.  
Consider a dense sequence $(\ba_n)_{n\in\N}\subset A^N$ (its existence is guaranteed by the fact that $A^N$ is a compact metric space), and define by recursion the sequence of measurable functions $\bpi_n:\Xc^N\rightarrow A^N$ as
\beqs 
\bpi_0(\bx)&=&\ba_0\\
\bpi_{n+1}(\bx)&=&\begin{cases}
\bpi_{n}(\bx) \text{ if } h(\bx,\bpi_{n}(\bx))\leq h(\bx,\ba_{n+1})\\
\ba_{n+1} \text{ else.}
\end{cases}
\enqs
The measurability of $\bpi_{n}$ is easily established by induction on $n$: For $n=0$, it is clear. Assuming that $\bpi_{n}$ is measurable, and denoting 
\beqs 
g_n(\bx) &= & h(\bx,\bpi_{n}(\bx))-h(\bx,\ba_{n+1}), \forall \bx\in\Xc^N,
\enqs 
notice that for any measurable set $B\subset A^N$, we have
\beqs 
[\bpi_{n+1}]^{-1}(B)=\begin{cases}
[\bpi_{n}]^{-1}(B)\cap g_n^{-1}(\R_-)  \text{ if }\ba_{n+1}\not\in B,\\
\big( [\bpi_{n}]^{-1}(B)\cap g_n^{-1}(\R_-)\big)\cup g_n^{-1}(\R_+^\star) \text{ if } \ba_{n+1}\in B,
\end{cases}
\enqs 
which is clearly a measurable set, and proves  the induction. Then, let us consider an embedding $\phi:A^N\rightarrow [0,1]$ such that $\phi$ and $\phi^{-1}$ are uniformly continuous (see Lemma C.2 in \cite{motte2019meanfield}). Then, $(\phi\circ \bpi_n)_{n\in\N}$ denotes a sequence of measurable functions from $\Xc^N$ to 
$\phi(A^N)\subset [0,1]$. It is well known that the function $\liminf_{n\in\N} (\phi\circ \bpi_n)$ is then measurable from $\Xc^N$ to $\phi(A^N)$ (we here use the fact that $\phi$ is continuous and $\phi(A^N)$ is closed, which ensures that the $\liminf$ takes its values in $\phi(A^N)$). Finally, let us denote $\bpi^\star:\Xc^N\rightarrow A^N$ defined by
\beqs 
\bpi^\star &=& \phi^{-1}\circ \liminf_{n\in\N} (\phi\circ \bpi_n).
\enqs 
$\bpi^\star$ is then measurable by composition. Furthermore, for any $\bx\in\Xc^N$, $\phi\circ\bpi^\star(\bx)=\liminf_{n\in\N} (\phi\circ \bpi_n(\bx))$ is an accumulation point of the sequence $(\phi\circ \bpi_n(\bx))_{n\in\N}$, which implies, by continuity of $\phi^{-1}$, that $\bpi^\star(\bx)$ is an accumulation point of $(\bpi_n(\bx))_{n\in\N}$. Given the definition of $\bpi_n(\bx)$, it is clear by induction that for any $n\in\N$, $h(\bx, \bpi_n(\bx))\leq \min_{m\leq n} h(\bx, \bpi^\star_m)$, and thus  $h(\bx, \bpi^\star(\bx))\leq \min_{n\in\N} h(\bx, \bpi_n)$. By density of $(\ba^n)_{n\in\N}$ and by continuity of $h(\bx,\cdot)$, this implies that 
$h(\bx, \bpi^\star(\bx))=\min_{\ba\in A^N} h(\bx, \ba)$ for all $\bx\in \Xc^N$, i.e. $\bpi^\star(\bx)\in \argmin_{\ba \in A^N} h(\bx, \ba)$. We conclude that $\bpi^\star$ is thus a measurable selection of $\argmin_{\ba \in A^N} h(., \ba)$. 
\ep

\vspace{3mm}

By Lemma \ref{remmesur}, there exists a randomized feedback policy $\bpi^{\mfa_\epsilon,N}$ s.t. 
\beqs
\Wc_\bomrd\big(\Lc(\xi_\bx,\mfa_\epsilon(\mu_{_N}[\bx],\xi_\bx,U)), \mu_{_N}[\bx,\bpi^{\mfa_\epsilon,N}] \big)  & = & \inf_{\ba \in A^N} \Wc_\bomrd\big(\Lc(\xi_\bx,\mfa_\epsilon(\mu_{_N}[\bx],\xi_\bx,U)), \mu_{_N}[\bx,\ba] \big),
%& \leq & 2 M_N, 
\enqs
and the r.h.s. of the above equality is bounded by $2M_N$ from  Lemma \ref{lem-action-couple}. Together with Proposition \ref{prop-policy-1}, this proves Theorem \ref{theo1-controlNMDP}.

\subsection{Proof of Theorem \ref{theo2-controlNMDP} }\label{sec-chaos-policy2}

Given a  randomized feedback policy  $\bpi_r$ $\in$  $L^0(\Xc^N\times[0,1]^N;A^N)$), the set of measurable functions from  $\Xc^N\times [0,1]^N$) into $A^N$, the associated feedback control is the unique control $\balpha^\pi$ given by  
$\balpha^{\bpi}_t$ $=$  $\bpi_r(\bX_t,\bU_t)$, $t\in\N$,   where $\{\bU_t$ $=$ $(U_t^i)_{i\in \llbracket 1,N\rrbracket}, t\in\N\}$ is a family of mutually i.i.d. uniform random variables on $[0,1]$, independent of $\Gc$, $\beps$. 
By misuse of notation, we denote $V^{\bpi}_N=V^{\balpha^{\bpi}}_N$. For $\bpi_r$ $\in$  $L^0(\Xc^N\times[0,1]^N;A^N)$, we introduce the operator $\Tc^{\bpi_r}_N$  on $L^\infty_m(\Xc^N)$, defined  by
\begin{align} 
\Tc_N^{\bpi_r}W(\bx) & : = \; \E[\bff(\bx, \bpi_r(\bx,\bU_0))+\beta W(\bF(\bx, \bpi_r(\bx,\bU_0), \beps_1)], \quad \forall \bx\in \Xc^N,
\end{align}
where $\bU_0$ $=$ $(U_0^i)_{i\in \llbracket 1,N\rrbracket}$ is a family of  i.i.d.  $\sim$ $\Uc([0,1])$, independent of $\Gc$, $\beps$.

\vspace{2mm}

We adapt  Proposition \ref{prop-policy-1} to  the case of   randomized feedback policies.

\begin{Proposition}\label{prop-policy-1-rand}
Let $\mfa_\epsilon$ be an $\eps$-optimal randomized feedback policy for the CMKV-MDP, and  consider any  feedback policy  $\bpi_r$ for the $N$-agent MDP.  Then, the feedback control $\balpha^{\bpi_r}$ is 
\beqs 
\Oc(\epsilon + \sup_{\bx\in\Xc^N} \E\Big[ \Wc_\bomrd\big(\Lc(\xi_\bx,\mfa_\epsilon(\mu_{_N}[\bx],\xi_\bx,U)), \mu_{_N}[\bx,\bpi_r(\bx,\bU_0)] \big) \Big]^\gamma + M_N^\gamma) \text{-optimal for } V_N(\bx_0),
\enqs 
namely
\beqs
V_N(\bx_0)  -  C \Big( \epsilon +   \sup_{\bx\in\Xc^N} \E\Big[ \Wc_\bomrd\big(\Lc(\xi_\bx,\mfa_\epsilon(\mu_{_N}[\bx],\xi_\bx,U)), \mu_{_N}[\bx,\bpi_r(\bx,\bU_0)] \big) \Big]^\gamma + M_N^\gamma \Big) & \leq &  V_N^{\bpi}(\bx_0). 
\enqs
Here $(\xi_\bx,U)$ $\sim$ $\mu_{_N}[\bx]\otimes\Uc([0,1])$,  and $\bU_0$ $=$ $(U_0^i)_{i\in\llbracket 1,N\rrbracket}$ is a family of i.i.d. $\sim$ $\Uc([0,1])$, independent of $\beps$.  
\end{Proposition}
{\bf Proof.}  Fix $\bx\in\Xc^N$, and let  $\ba$ $=$ $\bpi_r(\bx,\bU_0)$ be the random variable valued  in $A^N$. 
%and let $\mra_\eps$ $\in$ $L^0(\Xc\times [0,1];A)$ defined by $\mra_\eps(x,u)$ $=$ $\mfa_\eps(\mu_{_N}[\bx],x,u)$ for $x\in\Xc$, $u\in [0,1]$. 
By definition, we have  $\Tc_N^{\bpi_r}V_N(\bx)$ $=$ $\E\big[ \T^{\ba}_NV_N(\bx)\big]$.  By Lemma \ref{lem-verif-res}, we have 
\beqs 
\T^{\ba}_NV_N(\bx) &\geq& V_N(\bx)-\epsilon - C \big[   \Wc_\bomrd\big(\Lc(\xi_\bx,\mfa_\epsilon(\mu_{_N}[\bx],\xi_\bx,U)), \mu_{_N}[\bx,\bpi_r(\bx,\bU_0)] \big)^\gamma + M_N^\gamma \big].  
%\\ &\geq& V_N(\bx)-\eps-\Oc(\sup_{\bx\in\Xc^N}\Wc(\P_{\xi,a(\xi,U)}, \frac{1}{N}\sum_{i=1}^N \delta_{x^i,\ba^i})^\gamma + M_N^\gamma)
\enqs 
Taking the expectation, and by Jensen's inequality, we then get
\beqs 
\Tc_N^{\bpi_r} V_N(\bx)&\geq& V_N(\bx)-\epsilon - C \Big( \sup_{\bx\in\Xc^N}  \E \Big[ \Wc_\bomrd\big(\Lc(\xi_\bx,\mfa_\epsilon(\mu_{_N}[\bx],\xi_\bx,U)), \mu_{_N}[\bx,\bpi_r(\bx,\bU_0)] \big) \Big]^\gamma + M_N^\gamma \Big),
\enqs 
and we conclude by the verification result in Lemma \ref{lemverif-N}. 
\ep

 \vspace{3mm}

 Compared to Proposition \ref{prop-policy-1}, Proposition \ref{prop-policy-1-rand} means that with a randomized feedback policy $\bpi_r$, one can obtain a ``good" performance whenever it  produces  empirical state-action distributions that are close the theoretical state-action distribution generated by $\mfa_\epsilon$ {\it on average}, i.e., 
 that makes the quantity 
\beqs 
 \E \Big[ \Wc_\bomrd\big(\Lc(\xi_\bx,\mfa_\epsilon(\mu_{_N}[\bx],\xi_\bx,U)), \mu_{_N}[\bx,\bpi_r(\bx,\bU_0)] \big) \Big]
\enqs
as small as possible.  More precisely, if we can design a randomized policy $\bpi_r$ such that 
\beqs 
 \E \Big[ \Wc_\bomrd\big(\Lc(\xi_\bx,\mfa_\epsilon(\mu_{_N}[\bx],\xi_\bx,U)), \mu_{_N}[\bx,\bpi_r(\bx,\bU_0)] \big) \Big] & \leq &  C M_N, 
\enqs
then by Proposition \ref{prop-policy-1-rand}, this will prove the statement of Theorem \ref{theo2-controlNMDP}.  The next result shows how it can be achieved.

\begin{Lemma}\label{theo-policy}
%Let $K\in\R$ be a fixed constant. 
Let $\mfa$ $:$ $\Pc(\Xc)\times \Xc\times [0,1]\rightarrow A$ be any (if it exists) randomized feedback policy for the CMKV-MDP such that
\begin{align}\label{eq-cond-a}
\E[d_A(\mfa(\mu, x, U),  \mfa(\mu, x', U))]  & \leq \;   Kd(x,x'), \quad \forall  \mu \in \Pc(\Xc), x,x' \in \Xc, 
\end{align}
(here $U$ $\sim$ $\Uc[0,1])$) for some positive constant $K$. Consider the randomized feedback policy for the $N$-agent MDP defined by
\begin{align} \label{defoptpi} 
\bpi_r^{\mfa,N}(\bx,\bu) &= \;  \Big( \mfa(\mu_{_N}[\bx],x^i,u^i) \Big)_{i\in \llbracket 1,N\rrbracket}, \quad \bx = (x^i)_{i\in \llbracket 1,N\rrbracket} \in \Xc^N, \; \bu=(u^i)_{i\in \llbracket 1,N\rrbracket} \in [0,1]^N. 
\end{align}
Then, 
\beqs 
 \E \Big[ \Wc_\bomrd\big(\Lc(\xi_\bx,\mfa(\mu_{_N}[\bx],\xi_\bx,U)), \mu_{_N}[\bx,\bpi_r^{\mfa,N}(\bx,\bU_0)] \big) \Big] & \leq &  (2+K) M_N,  
\enqs
where $(\xi_\bx,U)$ $\sim$ $\mu_{_N}[\bx]\otimes\Uc([0,1])$.
\end{Lemma}
{\bf Proof.} Fix $\bx\in \Xc^N$, and  set $\mra_\bx(x,u)$ $=$ $\mfa(\mu_{_N}[\bx],x,u)$ for $(x,u)$ $\in$ $\Xc\times[0,1]$. 
Let us  consider a family $\bxi$ $=$ $(\xi^i)_{i\in\llbracket 1, N\rrbracket}$ of $N$ i.i.d. random variables such that $\xi^i$ $\sim$ $\mu_{_N}[\bx]$, and independent of  $\bU_0$. 
Let us consider $\sigma^{\bxi, \bx}$, the optimal permutation defined in Definition \ref{defpermut} between $\bxi$ and $\bx$. We have
\beqs 
& &  \Wc_\bomrd\big(\Lc(\xi_\bx,\mfa(\mu_{_N}[\bx],\xi_\bx,U)), \mu_{_N}[\bx,\bpi_r^{\mfa,N}(\bx,\bU_0)] \big)   \\
&= & \Wc_\bomrd\big(\Lc(\xi_\bx, \mra_\bx(\xi_\bx,U)),    \frac{1}{N} \sum_{i=1}^N \delta_{x^i,\mra_\bx(x^i,U_0^i)}  \big)\\
&\leq& \Wc_\bomrd\big(\Lc(\xi_\bx, \mra_\bx(\xi_\bx,U)), \frac{1}{N}\sum_{i=1}^N \delta_{\xi^{\sigma^{\bxi,\bx}_i}, \mra_\bx(\xi^{\sigma^{\bxi,\bx}_i}, U_0^i)} \big) 
+\Wc_\bomrd\big(\frac{1}{N}\sum_{i=1}^N \delta_{\xi^{\sigma^{\bxi,\bx}_i}, \mra_\bx(\xi^{\sigma^{\bxi,\bx}_i}, U_0^i)}, \frac{1}{N}\sum_{i=1}^N \delta_{x^i, \mra_\bx(x^i, U_0^i)} \big)\\
& \leq &  \Wc_\bomrd\big(\Lc(\xi_\bx, \mra_\bx(\xi_\bx,U)), \frac{1}{N}\sum_{i=1}^N \delta_{\xi^{i}, \mra_\bx(\xi^{i}, U_0^{(\sigma^{\bxi,\bx})^{-1}_i})} \big) 
+ \bomrd_N\big( (\bxi^{\sigma^{\bxi,\bx}}, \bpi_r^{\mfa,N}(\bxi^{\sigma^{\bxi,\bx}},\bU_0)),(\bx,\bpi_r^{\mfa,N}(\bx,\bU_0)) \big),     
%+  \frac{1}{N}\sum_{n=1}^N d(((\xi^{\sigma^{\bxi,\bx}_n)}, a(\xi^{\sigma^{\bxi,\bx}_n}, U^n)), (x^n, a(x^n, U^n)))
\enqs
where we set $\bxi^{\sigma^{\bxi,\bx}}$ $=$ $(\xi^{\sigma^{\bxi,\bx}_i})_{i\in\llbracket 1,N\rrbracket}$, and use \eqref{inegW} in the last inequality. 
Taking the expectation,  we then obtain  under condition \eqref{eq-cond-a}
\beqs 
& &  \Wc_\bomrd\big(\Lc(\xi_\bx,\mfa(\mu_{_N}[\bx],\xi_\bx,U)), \mu_{_N}[\bx,\bpi_r^{\mfa,N}(\bx,\bU_0)] \big)   \\
%&&\E[\Wc(\P_{\xi, a(\xi,U)}, \frac{1}{N}\sum_{n=1}^N \delta_{x^n, a(x^n, U^n)})]\\
%&\leq& \E[\Wc(\P_{\xi, a(\xi,U)}, \frac{1}{N}\sum_{n=1}^N \delta_{\xi^{n}, a(\xi^{n}, U^{(\sigma^{\bxi,\bx})^{-1}_n})})]+\frac{1}{N}\sum_{n=1}^N \E[d(((\xi^{\sigma^{\bxi,\bx}_n)}, a(\xi^{\sigma^{\bxi,\bx}_n}, U)), (x^n, a(x^n, U)))]\\
&\leq& M_N+(1+K) \E[ \bd_N(\bxi^{\sigma^{\bxi,\bx}},\bx)]  \\
& = &  M_N+(1+K) \E\big[\Wc_d\big(\mu_{_N}[\bxi], \mu_{_N}[\bx] \big) \big] \; \leq  \;  (2 + K) M_N,
\enqs 
where we use \eqref{couplingpermut} in the last equality.  This concludes the proof.
\ep

\vspace{3mm}

We now apply Lemma \ref{theo-policy} with an $\eps$-optimal randomized feedback policy $\mfa_\epsilon$ for the CMKV-MDP,  and combined with Proposition \ref{prop-policy-1-rand}, this proves the required result in Theorem  \ref{theo2-controlNMDP}.

\appendix

\section{Existence of optimal randomized control for CMKV-MDP} \label{sec-optcon}

Recall from Proposition 4.1 in \cite{motte2019meanfield} that the Bellman operator $\Tc$ of the CMKV-MDP  is written in the lifted form as 
\begin{align}  \label{defTc} 
[\Tc W](\mu) & = \; \sup_{\boa \in \bmrA} \Big\{  \tilde f(\mu,\boa) + \beta \E \big[ W\big( \tilde F(\mu,\boa,\eps_1^0)\big) \big] \Big\},  \quad \mu \in \Pc(\Xc), 
\end{align}
for $W$ $\in$ $L_m^\infty(\Pc(\Xc))$, where $\bmrA$ $=$ $\Pc(\Xc\times A)$,  $\tilde F$ is the measurable function on $\Pc(\Xc)\times\bmrA\times E^0$ $\rightarrow$ $\Pc(\Xc)$ defined by 
\begin{align}
\tilde  F(\mu,\boa,e^0) &= \; F(\cdot, \cdot,\bp(\mu,\boa),\cdot,e^0)\star\big(\bp(\mu,\boa)\otimes \Lc(\eps_1) \big), 
\end{align}
and $\tilde f$ is the measurable function on $\Pc(\Xc)\times\bmrA$ defined by
\begin{align}
\tilde f(\mu,\boa) & = \;  \int_{\Xc\times A}  f(x,a,\bp(\mu,\boa))\bp(\mu,\boa)(\d x,\d a).   
\end{align}
Here  $\star$ is the pushforward measure notation, $\bp$ is a measurable coupling projection from $\Pc(\Xc)\times\bmrA$ into $\bmrA$: $\bp(\mu,\boa)=\bp(\mu,\bp(\mu,\boa))$, satisfying 
$\text{pr}_{_1}\star \bp(\mu,\boa) \; = \; \mu$, and  $\bp(\mu,\boa)=\boa$ if $\text{pr}_1\star \boa=\mu$ (where $\text{pr}_1$ is the projection function on the first coordinate).  
Since $\tilde f$ and $\tilde F$ depend upon $\boa$ only through $\bp(\mu,\boa)$, it is clear that the supremum  in \eqref{defTc}, for each $\mu$ $\in$ $\Pc(\Xc)$,  can be taken actually over the the subset 
$\Gamma_\mu$ $:=$ $\{\boa: (\mu,\boa)\in\Gamma\}$ $\subset$ $\bmrA$, where $\Gamma:=\{(\mu, \boa)\in \Pc(\Xc)\times \bmrA: \text{pr}_1 \star \boa = \mu\}$ is closed in $\Pc(\Xc)\times \bmrA$ from the continuity of  $\boa$ $\mapsto$ $\text{pr}_1\star \boa$. 
%as the graph of the continuous function $\text{pr}_1$. 
Moreover, since $V$ is  continuous (see \eqref{Vholder}), it is straightforward to prove that 
\beqs 
(\mu,\boa) \in \Gamma &\mapsto&  \tilde f(\mu,\boa) + \beta \E \Big[ V\big( \tilde F(\mu,\boa,\eps_1^0)\big) \Big] \\
& & \; = \; \int_{\Xc\times A} f(x,a,\boa) \boa(\d x,\d a)  + \beta \E \Big[ V\big( F(\cdot, \cdot,\mu,\cdot,e^0)\star\big(\boa\otimes \Lc(\eps_1) \big) \big) \Big] 
\enqs 
is continuous and thus upper continuous on $\Gamma$.  Therefore,  by \cite{bertsekas1996stochastic}, Proposition 7.33, there exists a measurable function $\phi:\Pc(\Xc)\rightarrow \bmrA$ whose graph is included in  $\Gamma$ and such that
\begin{align} 
\tilde f(\mu,\phi(\mu)) + \beta \E \big[ V\big( \tilde F(\mu,\phi(\mu),\eps_1^0)\big) \big] &= \;  \sup_{\boa \in \Gamma_\mu} \Big\{  \tilde f(\mu,\boa) + \beta \E \big[ V\big( \tilde F(\mu,\boa,\eps_1^0)\big) \big] \Big\},  \nonumber \\
&= \;  [\Tc V](\mu) \; = V(\mu), \quad \forall \mu\in\Pc(\Xc), \label{Vsup} 
\end{align}  
where the last equality follows from the fixed point equation of $V$.  By the universal disintegration theorem (see \cite{kallenberg2017random}, Corollary 1.26), there exists $\kappa:\Xc\times \Pc(\Xc\times A)\times \Pc(\Xc)\rightarrow \Pc(A)$ such that for all $\boa\in\Pc(\Xc\times A)$, $\mu\in\Pc(\Xc)$ with 
$\text{pr}_1\star \boa$ $=$ $\mu$, we have $\boa=\mu\hat{\otimes}\kappa(\cdot, \boa,\mu)$ (where $\hat{\otimes}$ denotes the probability-kernel product). Furthermore, by Blackwell-Dubins Lemma,  there exists a measurable function $\rho:\Pc(A)\times [0,1]\rightarrow A$ such that for all $\pi\in\Pc(A)$, if $U$ denotes a uniform random variable, then 
$\rho(\pi,U)\sim \pi$. We can then define the randomized feedback policy
\beqs 
\mfa_0(\mu,x,u) &=& \rho(\kappa(x,\phi(\mu), \mu), u), 
\enqs 
which satisfies by construction $\Lc(\xi,\mfa_0(\mu,\xi,U))$ $=$ $\phi(\mu)$ for $(\xi,U)$ $\sim$  $\mu\otimes\Uc([0,1])$ so that 
\beqs
\tilde f(\mu,\phi(\mu)) &=& \E \Big[ f\big(\xi,\mfa_0(\mu,\xi,U), \Lc(\xi,\mfa_0(\mu,\xi,U)) \big) \Big]  \\
\tilde F(\mu,\phi(\mu),\eps_1^0) & = & \P^0_{F(\xi,\mfa_0(\mu,\xi,U), \Lc(\xi,\mfa_0(\mu,\xi,U)),\eps_1,\eps_1^0)}. 
\enqs
Recalling notation in \eqref{defTa}, and by \eqref{Vsup}, this shows that 
\beqs
\T^{\mfa_0(\mu,.)} V(\mu)=V(\mu).
\enqs
According to the verification result (Proposition 4.3 in \cite{motte2019meanfield}), this ensures that that the randomized feedback control $\alpha^0$ $\in$ $\Ac$ defined by 
\begin{align} \label{ranfeedback0} 
\alpha_t^0 &= \;  \mfa_0(\P_{X_t}^0,X_t,U_t), \quad t \in \N, 
\end{align} 
where $(U_t)_{t\in\N}$ is an i.i.d. sequence of random variables, $U_t$ $\sim$ $\Uc([0,1])$, independent of $\xi_0$ $\sim$ $\mu_0$, and $\beps$,  is  an  optimal control for $V(\mu_0)$.

\section{Bellman equation for the $N$-agent MDP}\label{sec-Bell-N}

In this section, we study and rigorously state  properties on the Bellman equation for the $N$-agent problem,  viewed  as a MDP  with state space $\Xc^N$, action space $A^N$, noise sequence $\beps=(\beps_t)_{t\in\N^\star}$ with  $\beps_t$ $:=$ $((\eps^i_t)_{i\in\llbracket 1, N\rrbracket}, \eps^0_t)$  valued in $E^N\times E^0$, 
state transition function
\beqs 
\bF(\bx, \ba, \be) &:= & \Big(F(x^i, a^i, \mu_{_N}[\bx,\ba], e^i, e^0)\Big)_{i\in \llbracket 1,N\rrbracket},  \quad  \be=((e^i)_{i\in \llbracket 1,N\rrbracket}, e^0)\in E^N\times E^0,
\enqs 
and reward function 
\beqs
\bff (\bx,\ba) &= &\frac{1}{N}\sum_{i=1}^N f\big(x^i, a^i, \mu_{_N}[\bx,\ba] \big), \quad \bx = (x^i)_{i\in \llbracket 1,N\rrbracket}, \; \ba = (a^i)_{i\in \llbracket 1,N\rrbracket}. 
\enqs 
With respect to standard framework of MDP,  we pay a careful attention when dealing with possibly continuous state/action spaces $(\Xc,A)$,   and optimizing in general over open-loop controls.

Let us consider  the set $\bVc$ of sequences $\bnu=(\bnu_t)_{t\in\N}$ with $\bnu_0$ a measurable function from $([0,1]^{N})^{\N}$ into $A^N$, and 
$\bnu_t$ a measurable function from $([0,1]^{N})^{\N}\times (E^N\times E^0)^t$ into $A^N$ for $t$ $\in$ $\N^\star$. For each $\bnu$ $\in$ $\bVc$, we can associate a control process $\balpha^{\bnu}$ $\in$ $\bAc$ given by
 \beqs 
\balpha^{\bnu}_t & := &  \bnu_t(\bU, (\beps_s)_{s\in\llbracket 1, t\rrbracket}), \quad  t\in\N,
\enqs 
(with the convention that $\balpha_0^{\bnu}$ $=$ $\bnu_0(\bU)$ when $t$ $=$ $0$), where $\bU$ $=$ $(U_t^i)_{i\in\llbracket 1,N\rrbracket, t \in \N}$ is a family of mutually i.i.d. uniform random variables on $[0,1]$, independent of $\beps$, and conversely any control $\balpha \in \bAc$ can be represented as $\balpha^{\bnu}$ for some $\bnu\in \bVc$. 
We call $\bVc$ the set of randomized open-loop policies. By misuse of notation, we write  $V_N^{\bnu}$ $=$ $V_N^{\balpha^{\bnu}}$.

\vspace{1mm}

Let us denote by $L^\infty(\Xc^N)$ the set of bounded real-valued functions on $\Xc^N$, 
and by $L^\infty_m(\Xc^N)$ the subset of measurable functions in $L^\infty(\Xc^N)$. 
We then introduce the Bellman ``operator'' $\Tc_N:L^\infty_m(\Xc^N)\rightarrow L^\infty(\Xc^N)$ defined 
for any $W$ $\in$ $L^\infty_m(\Xc^N)$ by: 
\begin{align} \label{defTc-N} 
[\Tc_N W](\bx) & := \; \sup_{\ba \in A^N} \T^{\ba}_N W (\bx),  \quad \bx \in \Xc^N. \end{align}
where
\begin{align} \label{defTT-N} 
\T^{\ba}_N W (\bx) & := \; \bff(\bx,\ba) + \beta \E \big[ W\big(  \bF(\bx,\ba, \beps_1)\big) \big],  \quad \bx \in \Xc^N,  \; \ba \in A^N. 
\end{align}
Notice that the $\sup$ can a priori lead to a non measurable function $\Tc_N W$. Because of this, $\Tc_N$ is not an operator on $L^\infty_m(\Xc^N)$ in the strict sense. To see $\Tc_N$ as an operator, we have to find a subset in   $L^\infty_m(\Xc^N)$ that is preserved by $\Tc_N$. The next result introduces such subset.

\begin{Lemma}\label{lem-Tc-reg}\label{lem-Tc-M}
Let $\Mc$ be the set in $L^\infty_m(\Xc^N)$ defined by 
\begin{align}\label{eq-M}
\Mc &:= \;  \Big\{  W\in L^\infty_m(\Xc^N): \big\vert W(\bx) - W(\bx') \big\vert \leq  2K_f \sum_{t=0}^\infty \beta^t \min \big[ (2K_F)^t \bd_N (\bx,\bx'),\Delta_\Xc \big], \;  \forall  \bx,\bx' \in \Xc^N\}.
\end{align} 
Then $\Mc$ is a complete metric space under the $\Vert\cdot \Vert$ norm, and $\T^{\ba}_N$, for all $\ba\in A^N$, and $\Tc_N$, preserve $\Mc$: 
$\T^{\ba}_N \Mc$ $\subset$ $\Mc$, $\Tc_N \Mc$ $\subset$ $\Mc$. 
\end{Lemma}
\noindent {\bf Proof.}
It is clear that $\Mc$ is closed in $L^\infty_m(\Xc^N)$, and is therefore a complete metric space for $\Vert\cdot \Vert$. Let $W\in \Mc$. Fix $\bx,\bx'\in \Xc^N$, and  $\ba\in A^N$. Let us start with two preliminary estimations: under (${\bf Hf_{lip}}$), and recalling \eqref{inegW},  we clearly have
\begin{align}
\vert \bff(\bx, \ba)-\bff(\bx', \ba)\vert  &\leq  2 K_fd_N(\bx,\bx'). \label{f-estim-N}
\end{align}
Similarly, under (${\bf HF_{lip}}$), for $e^0\in E^0$, we have
\begin{align}\label{F-estim-N}
\E [\bd_N(\bF(\bx,\ba,(\eps^i_1)_{i\in \llbracket 1, N\rrbracket},e^0),\bF(\bx',\ba,(\eps^i_1)_{i\in \llbracket 1, N\rrbracket},e^0))]&\leq 2 K_F \bd_N(\bx,\bx').
\end{align}
Thus, denoting by  $\bX_1=\bF(\bx,\ba,(\eps^i_1)_{i\in \llbracket 1, N\rrbracket},e^0)$ and $\bX'_1=\bF(\bx',\ba,(\eps^i_1)_{i\in \llbracket 1, N\rrbracket},e^0)$, we have, by Jensen's inequality and then \eqref{F-estim-N},
\begin{align}
    \E\left[\sum_{t=0}^\infty \beta^t \min\big[ (2K_F)^t  \bd_N (\bX_1,\bX'_1),\Delta_\Xc \big] \right]&\leq \;  \sum_{t=0}^\infty  \beta^t \min\big[ (2K_F)^t \E[\bd_N (\bX_1,\bX'_1)],\Delta_\Xc \big] \nonumber \\
%&\leq \; \sum_{t=0}^\infty \beta^t \min\big[ (2K_F)^t 2K_F d_N(\bx,\bx'),\Delta_\Xc \big]  \nonumber \\
&\leq \;  \sum_{t=0}^\infty \beta^t \min\big[ (2K_F)^{t+1} \bd_N(\bx,\bx'),\Delta_\Xc\big].  \label{sum_estim-N} 
\end{align}
The definition of $\T^{\ba}_N$ combined with \eqref{f-estim-N}, the fact that $W\in\Mc$, and \eqref{sum_estim-N}, implies that
\beqs
\vert \T_N^{\ba} W(x)  -  \T_N^{\ba} W(x')\vert  &\leq & 2K_f  \bd_N(\bx,\bx')+\beta  2K_f\sum_{t=0}^\infty  \beta^t \min\big[ (2K_F)^{t+1} \bd_N(\bx,\bx'),\Delta_\Xc\big] \\
& \leq & 2 K_f  \sum_{t=0}^\infty  \beta^t \min\big[ (2K_F)^{t} \bd_N(\bx,\bx'),\Delta_\Xc\big],
\enqs
which shows that  $\T^{\ba}_N W\in \Mc$, i.e. $\T^{\ba}_N$ preserves $\Mc$. Furthermore, we have
\beqs 
\vert \Tc_N W(\bx) - \Tc_N W(\bx')\vert 
%&=&\vert \sup_{\ba \in A^N} \T_N^{\ba} W(\bx)- \sup_{\ba \in A^N} \T_N^{\ba} W(\bx')\vert\\
&\leq& \sup_{\ba \in A^N} \vert \T_N^{\ba} W(\bx)-\T_N^{\ba} W(\bx')\vert\\
%&\leq& \sup_{\ba \in \bA} \big( 2K_f \sum_{t=0}^\infty \beta^t \min((2K_F)^t d_N (\bx,\bx'),\Delta_\Xc)\big)\\
& \leq & 2K_f \sum_{t=0}^\infty \beta^t \min\big[ (2K_F)^t \bd_N (\bx,\bx'),\Delta_\Xc\big],
\enqs 
which also shows that  $\Tc_N W\in \Mc$.
\ep

\vspace{3mm}

Lemma \ref{lem-Tc-M} implies that by restricting $\Tc_N$ and $\T^\ba_N$ to $\Mc$, we can see $\Tc_N$ and $\T^\ba_N$ as operators on $\Mc$, that is, $\Tc_N:\Mc\rightarrow \Mc$ and $\T^\ba_N:\Mc\rightarrow \Mc$. However, the property defining the functions in $\Mc$ (see \eqref{eq-M}) is not very natural and practical. The following result provides a more convenient property satisfied by all functions in $\Mc$.

\begin{Lemma}\label{lem-M-hol}
%Let $\gamma:=\min\big[1,\frac{\vert \ln \beta\vert}{\ln (2 K_F)}\big]$. 
There exists $K_\star\in\R$ such that any function $W\in\Mc$ is $\gamma$-Hölder with constant factor $K_\star$, i.e.
\beqs 
\big\vert W(\bx) - W(\bx') \big\vert  &\leq&  K_\star \bd_N(\bx,\bx')^\gamma, \quad \forall  \bx,\bx' \in \Xc^N.
\enqs 
\end{Lemma}
\noindent {\bf Proof.}
We have
\beqs
\vert W(\bx)-W(\bx')\vert  & \leq &  2K_f \sum_{t=0}^\infty \beta^t \min\big[ (2K_F)^t \bd_N (\bx,\bx'),\Delta_\Xc\big] \; =: \; 2 K_f  S(\bd_N (\bx,\bx')).
\enqs
where $S(m)$ $=$ $\sum_{t=0}^\infty \beta^t \min[ (2K_F)^t m,\Delta_\Xc]$. If $2\beta K_F<1$, we clearly have
\beqs
S(m) &\leq&   m\sum_{t=0}^\infty (\beta 2K_F)^t \; = \;  \frac{m}{1-\beta 2K_F},
\enqs
and so $W$ is $1$-Hölder. Let us now study the case $2\beta K_F>1$. In this case, in particular, $2K_F>1$ since $\beta$ $\in$ $(0,1)$, thus  
$t$ $\mapsto$ $(2K_F)^t$ is nondecreasing, and so 
\begin{align}
S(m) &\leq\; \sum_{t=0}^\infty \int_t^{t+1}\beta^t   \min \big[ (2K_F)^s m, \Delta_\Xc \big]ds \\ 
&\leq\; \frac{1}{\beta}\sum_{t=0}^\infty \int_t^{t+1}\beta^s   \min \big[  (2K_F)^s m , \Delta_\Xc \big]ds  \; =  \; \frac{1}{\beta} \int_0^{\infty}e^{ s  \ln\beta }   \min \big[ m e^{s \ln(2K_F)}, \Delta_\Xc \big]ds. \label{serieint-N}
\end{align}
Let $t_\star$ $=$ $t_\star(m)$  be such that $me^{t_\star \ln (2K_F)}=\Delta_\Xc$, i.e.
%\beqs
$t_\star$ $=$ $\frac{\ln(\Delta_\Xc/m)}{\ln(2 K_F)}$. 
%\enqs
Then, 
\beqs
\int_0^{\infty}e^{s  \ln\beta }   \min \big[ m e^{s \ln(2K_F)}, \Delta_\Xc \big]ds &\leq& m\int_0^{t^\star} e^{s \ln (2K_F \beta)} ds+\Delta_\Xc\int_{t^\star}^\infty e^{s \ln (\beta) }  ds\\
& = & \frac{m}{\ln (2K_F \beta)}\left[ e^{t_\star \ln (2K_F \beta) }-1 \right] - \frac{\Delta_\Xc}{\ln \beta}e^{\ln (\beta)  t^\star} \\
& = &  \frac{m}{\ln (2K_F \beta)}\Big[ \Big(\frac{\Delta_\Xc}{m}\Big)^{\frac{\ln (2K_F\beta)}{\ln(2 K_F)}} -1 \Big]-\frac{\Delta_\Xc}{\ln \beta}\left(\frac{\Delta_\Xc}{m}\right)^{ \frac{\ln (\beta)}{\ln(2 K_F)}} \\
& = &  \Delta_\Xc\Big(\frac{1}{\ln (2K_F \beta)}-\frac{1}{\ln \beta} \Big)\Big(\frac{\Delta_\Xc}{m}\Big)^{ \frac{\ln (\beta)}{\ln(2 K_F)}} -\frac{m}{\ln (2K_F \beta)} \\
&  \leq  & C m^{\min\big[1,\frac{\vert \ln \beta\vert}{\ln (2 K_F)}\big]} = C m^{\gamma}, 
\enqs
for some positive constant $C$ depending on $K_F$, $\beta$ and $\Delta_\Xc$. This implies that $W$ is $\gamma$-Hölder with a constant factor $K_\star$ that is clearly independent of $W\in\Sc$. This concludes the proof.
\ep

\vspace{3mm}

The consequence of Lemmas \ref{lem-Tc-M} and \ref{lem-M-hol} is that the set $\Mc\subset L^\infty_m(\Xc)$ is a closed set, preserved by $\Tc_N$ and contains only functions that are $\gamma$-Hölder with factor $K_\star$. We are now  able to get the existence of a unique fixed point to the Bellman operator  $\Tc_N$.

\begin{Proposition}\label{lemTc-N}\label{proregul-N} 
%Assume that  $({\bf H_{lip}})$ holds true.
(i)  The operator $\Tc_N$ is monotone increasing: for $W_1,W_2\in L_m^\infty(\Xc^N)$, if $W_1$ $\leq$ $W_2$,  then $\Tc_N W_1$ $\leq$ $\Tc_N W_2$.
(ii) Furthermore, it is contracting on $L^\infty_m(\Xc^N)$ with Lipschitz factor $\beta$, and admits a unique fixed point 
in $L^\infty_m(\Xc^N)$, denoted by $V^\star_N$, hence solution to: 
\begin{align}
V^\star_N & = \; \Tc_N V^\star_N. 
\end{align}
Moreover,  $V^\star_N\in\Mc$, and thus $V^\star_N$ is $\gamma$-H\"older with constant factor $K_\star$.
%{\color{blue}Should we explicit the Holder multiplicative constant? Because from the %proof of theorem 2.1 it seems that the constant is a bit complicated, and would make the %result less elegant... Maybe we should just say that it only depends upon $K_F$, %$\beta$, and $\Delta_\Xc$, and that an explicit expression can be retrieved from the %proof?}
%{\green pas besoin d'expliciter cette constante}
%with a factor $\frac{2K}{1-(2K)^\gamma\beta}$.
%\noindent (ii) For all $\epsilon>0$, there exists $\eta>0$ s.t. if $x,x'\in\Xc^N$ with $d_N(x,x')<\eta$, then (\red{reprendre l'enonce}). 
%\begin{align} 
%d_N\big(\tilde{F}(x,\ba,e^0), \tilde{F}(x',\ba,e^0) \big) + \big|  \tilde{f}(x,\ba) -   \tilde{f}(x',\ba) \big| &  \leq \;  \epsilon, \quad \forall \ba \in \bA, e^0 \in E^0.  
%\end{align}
\end{Proposition}
\noindent {\bf Proof.}  (i) The monotonicity of $\Tc_N$ is clear. (ii)  The  $\beta$-contraction property of $\Tc_N$ is  obtained by standard arguments, which implies the uniqueness of a fixed point (but not the existence). Let us prove the existence of a fixed point. As $\Mc$ is preserved by $\Tc_N$, and is closed for $\Vert \cdot\Vert$, and therefore complete (as a closed subset of the complete space $L^\infty_m(\Xc^N)$), by the Banach fixed point theorem, $\Tc_N$ admits a unique fixed point $V^\star_N$ in $\Mc$.  By Lemma \ref{lem-M-hol},  this implies that $V^\star_N$ is $\gamma$-Hölder with constant factor $K_\star$, and concludes the proof.
\ep

\begin{Remark} 
{\rm  Notice that the above arguments  would not work if we considered, instead of $\Mc$, directly the set of $\gamma$-Hölder continuous functions. Indeed, while it is true that such set is stabilized by $\Tc_N$ (it essentially follows from \eqref{f-estim-N} and \eqref{F-estim-N}), the set of $\gamma$-Hölder continuous functions is not closed in $L^\infty_m(\Xc^N)$ (and thus not a complete metric space): there might indeed exist a converging sequence of $\gamma$-Hölder continuous functions with multiplicative factors (in the Hölder property) tending toward infinity, such that the limit function is not $\gamma$-Hölder anymore.
}
\ep
\end{Remark}

\vspace{3mm}

As a consequence of  Proposition \ref{lemTc-N}, we can  show the following relation between the value function $V_N$ of the $N$-agent MDP, and the fixed point $V^\star_N$ of the Bellman operator $\Tc_N$.

\begin{Lemma} \label{lemVstar-N}
 For all $\bx$ $\in$ $\Xc^N$,
 %and $\xi$ $\in$ $L^0(\Gc;\Xc)$  with $\Lc(\xi)$ $=$ $x$, 
 we have $V_N(\bx)\leq V^\star_N(\bx)$.
 \end{Lemma}
 \noindent{\bf Proof.}
 For  any $\bx$ $\in$ $\Xc^N$, $\bnu\in \bVc$, we have
 \begin{align}  
 \E \Big[\bff(\bx,\bnu_0(\bU))  + \beta  V^{\star}_N\big(\bF (x, \bnu_0(\bU), \beps_1) \big) \Big]&= \; \E \Big[\Big\{\bff(\bx,\bnu_0(u))  + \beta\E[  V^{\star}_N(\bF (\bx, \bnu_0(\bu), \beps_1))]\Big\}_{\bu:=\bU}\Big] \nonumber \\
 &= \; \E \Big[\T^{\bnu_0(\bU)}V^\star_N(\bx)\Big]\leq \Tc V^\star_N(\bx)=V^\star_N(\bx). \label{Vinter} 
 \end{align} 
 For any $(\bu, \be)\in ([0,1]^{N})^{\N}\times (E^N\times E^0)$, and for any $\bnu\in \bVc$, we define $\vec{\bnu}^{\bu,\be}\in \bVc$ by
\beqs 
\vec{\bnu}^{\bu,\be}_t(\bu', (\be'_s)_{s\in \llbracket 1, t\rrbracket}):= \bnu_{t+1}(\bu, \be, (\be'_s)_{s\in \llbracket 1, t\rrbracket}), \quad  (\bu', (\be'_s)_{s\in \llbracket 1, t\rrbracket})\in  ([0,1]^{N})^{\N}\times (E^N\times E^0)^t,  t\in\N.
\enqs 
Standard Markov arguments imply the following flow property for randomized open-loop policies:
 \beqs 
 V^{\bnu}_N(\bx) &=&  \E\Big[\bff(\bx,\bnu_0(\bU))  + \beta  V^{\vec{\bnu}^{\bU, \beps_1}}_N(\bF (\bx, \bnu_0(\bU), \beps))\Big].
 \enqs 
Together with \eqref{Vinter}, we then get
 \beqs 
 V^\star_N(\bx)-V^{\bnu}_N(\bx)&\geq& \beta\E \Big[  V^{\star}_N(\bF (\bx, \bnu_0(\bU), \beps_1)-V^{\vec{\bnu}^{\bU, \beps_1}}_N(\bF (\bx, \bnu_0(\bU), \beps_1))\Big]\\
 &\geq& \beta \inf_{\bx \in \Xc^N,\bnu\in \bVc} \big\{ V^\star_N(\bx)-V^{\bnu}_N(\bx) \big\}.
 \enqs 
 Taking the infimum over $\bx \in \Xc^N, \bnu\in\bVc$ on the left hand side of the above inequality, and since $\beta$ $<$ $1$, this shows that 
 %\beqs 
 %\inf_{\bx \in \Xc^N,\bnu\in \bVc} \big\{ V^\star_N(\bx)-V^{\bnu}_N(\bx) \big\}\geq 0
 %\enqs
 %and thus 
 $V^{\bnu}_N(\bx)\leq V^\star_N(\bx)$ for all $\bnu\in\bVc$. We conclude that $V_N$ $\leq$ $V^\star_N$. 
\ep

\vspace{3mm}

We  aim now  to prove rigorously  the equality $V_N$ $=$ $V^\star_N$, i.e., the value function $V_N$ of the $N$-agent MDP satisfies the Bellman fixed point equation: $V_N$ $=$ $\Tc_N V_N$, 
and also to show the existence of $\eps$-optimal (randomized)  feedback control for $V_N$.

A  feedback policy (resp. randomized feedback policy)  is an element $\bpi$ $\in$ $L^0(\Xc^N;A^N)$ (resp. $L^0(\Xc^N\times[0,1]^N;A^N)$), the set of measurable functions from $\Xc^N$ (resp. $\Xc^N\times [0,1]^N$) into $A^N$. The associated feedback control is the unique control $\balpha^\pi$ given by  
$\balpha^{\bpi}_t$ $=$ $\bpi(\bX_t)$,  (resp. $\bpi_r(\bX_t,\bU_t)$), $t\in\N$,   where $\{\bU_t$ $=$ $(U_t^i)_{i\in \llbracket 1,N\rrbracket}, t\in\N\}$ is a family of mutually i.i.d. uniform random variables on $[0,1]$, independent of $\Gc$, $\beps$. 
By misuse of notation, we denote $V^{\bpi}_N=V^{\balpha^{\bpi}}_N$. Given $\bpi$ $\in$  $L^0(\Xc^N;A^N)$ (resp. $L^0(\Xc^N\times[0,1]^N;A^N)$), we introduce the operator $\Tc^{\bpi}_N$  on $L^\infty_m(\Xc^N)$, defined  by
\begin{align} \label{defTpi-N} 
\Tc^{\bpi}_N W(\bx) & := \;  \bff( \bx,\bpi(\bx)) + \beta \E \big[ W\big( \bF(\bx,\bpi(\bx),\beps_1)\big) \big], \quad \bx \in \Xc^N,  
\end{align}
resp. 
\begin{align}\label{defTpiR-N}
\Tc_N^{\bpi}W(\bx) & : = \; \E[\bff(\bx, \bpi(\bx,\bU_0))+\beta W(\bF(\bx, \bpi(\bx,\bU_0), \beps_1)], \quad \forall \bx\in \Xc^N,
\end{align}
where $\bU_0$ $=$ $(U_0^i)_{i\in \llbracket 1,N\rrbracket}$ is a family of  i.i.d.  $\sim$ $\Uc([0,1])$, independent of $\Gc$, $\beps$.

\vspace{3mm}

We have the basic and standard  properties on the operator $\Tc_N^{\bpi}$:

\begin{Lemma} \label{lemTcpi-N} 
Fix $\bpi$ $\in$ $L^0(\Xc^N;A^N)$ (resp. $L^0(\Xc^N\times[0,1]^N;A^N)$). 
\begin{itemize}
\item[(i)] The operator $\Tc^{\bpi}_N$ is $\beta$-contracting on $L^\infty_m(\Xc^N)$, and $V_N^{\bpi}$ is its unique fixed point.  
\item[(ii)] Furthermore, it is monotone increasing: for $W_1,W_2\in L^\infty(\Xc^N)$, if $W_1$ $\leq$ $W_2$,  then $\Tc^{\bpi}_N W_1$ $\leq$ $\Tc^{\bpi}_N W_2$.
\end{itemize}
\end{Lemma}

We state the standard verification type result for the $N$-individual MDP, by means of the Bellman operator.

\begin{Lemma}\label{lemverif-N} 
{\bf (Verification result)} 

\noindent Fix $\epsilon\geq0$, and suppose that  there exists an $\epsilon$-optimal  (randomized) feedback policy $\bpi^\epsilon$ 
%lifted randomized feedback policy  $\bpi_\epsilon$ $=$ $\bpi^{\mfa_\epsilon}$, $\mfa_\eps$ $\in$ $L^0(\Xc^N\times \Xc\times [0,1];A)$, 
for $V^\star_N$ in the sense that
\begin{align}
V^\star_N &\leq \;  \Tc^{\bpi^\epsilon}_N V^\star_N  +\epsilon.   
\end{align}
Then, $\balpha^{\pi^\epsilon}$ $\in$ $\bAc$  is $\frac{\epsilon}{1-\beta}$-optimal for $V_N$, i.e.,  $V^{\bpi^\epsilon}_N$ $\geq$ $V_N - \frac{\epsilon}{1-\beta}$, 
and we have $V_N$ $\geq$ $V^\star_N - \frac{\epsilon}{1-\beta}$. 
\end{Lemma}
{\bf Proof.} Since $V^{\bpi^\epsilon}_N$ $=$ $\Tc^{\bpi^\epsilon}_N V^{\bpi^\epsilon}_N$, and recalling from Lemma \ref{lemVstar-N} that $V^\star_N$ $\geq$ $V_N$ $\geq$ 
$V^{\bpi^\epsilon}_N$, we have for all $\bx$ $\in$ $\Xc^N$,
\begin{align}
 \Big| (V^\star_N - V^{\bpi^\epsilon}_N)(\bx) \Big|  &\leq \; \Big|  \Tc^{\bpi^\epsilon}_N V^\star_N(\bx)-\Tc^{\bpi^\epsilon}_N V^{\bpi^\epsilon}_N(\bx)\Big|+\epsilon  \; \leq \;  \beta \Vert V^\star_N-V^{\bpi^\epsilon}_N\Vert +\epsilon,
\end{align} 
where we used the $\beta$-contraction property of $\Tc_N^{\bpi_\epsilon}$ in Lemma \ref{lemTcpi-N}.   We deduce that  
$\Vert V^\star_N -V^{\bpi^\epsilon}_N \Vert$ $\leq$ $\frac{\epsilon}{1-\beta}$,  and then, $V_N$ $\geq$  $V^{\bpi^\epsilon}_N$ $\geq$ 
$V^\star_N-\frac{\epsilon}{1-\beta}$, which combined with $ V^\star_N\geq V_N$, concludes the proof.  
\ep

%\begin{Corollary}
%Notice that if for all $\eps>0$ we can find such "feedback policy" $a$, it would prove that in the MDP on $\Xc^N$, we could restrict ourselves to feedback controls, but it would also prove the DPP %\eqref{DPP} since we would have $\hat{V}=V^\star$. 
%\end{Corollary}

\vspace{3mm}

We finally  conclude this section  by showing the existence of an $\eps$-optimal (randomized)  feedback policy for $N$-agent MDP on $\Xc^N$, and obtain as a by-product the corresponding Bellman fixed point equation for its value function.

\begin{Proposition} \label{theomainDPP-N} 
%Assume that  $({\bf H_{lip}})$ holds true.   
For all $\epsilon>0$,  there exists a (randomized) feedback policy $\bpi^\epsilon$ that is 
$\epsilon$-optimal for $V^\star_N$.   Consequently, the control $\balpha^{\bpi^\epsilon}$ $\in$ $\bAc$ is $\frac{\epsilon}{1-\beta}$-optimal for $V_N$, and we have 
$V_N$ $=$ $V^\star_N$, which thus satisfies the Bellman fixed point equation.   
\end{Proposition}
{\bf Proof.}  We prove the result for $\epsilon$-optimal  feedback policy (the case of $\epsilon$-optimal randomized feedback policy is dealt with similarly). 
Fix $\epsilon$ $>$ $0$, and given $\eta$ $>$ $0$, consider a quantizing grid $\Mc^\eta$ $=$ $\{\bx_1,\ldots,\bx_{N^\eta}\}$ $\subset$ $\Xc^N$, and an associated partition  $C^\eta_k$, $k$ $=$ $1,\ldots,N^\eta$, of $\Xc^N$, satisfying 
\begin{align}
C^\eta_k & \subset \; B^\eta(\bx_k) := \Big\{ \bx \in \Xc^N: \bd_N(\bx,\bx_k) \leq \eta \Big\}, \quad k=1,\ldots,N_\eta. 
\end{align}
For any $\bx_k$, $k$ $=$ $1,\ldots,N^\eta$,  
%by \eqref{VepsY},
there exists $\ba^\eps_k$ $\in$ $A^N$ such that 
\begin{align} \label{VepsYi-N} 
V^\star_N(\bx_k) & \leq \T^{\ba^\eps_k} V_N^\star(\bx_k) + \frac{\epsilon}{3}.  
%\E \Big[ f(Y^{i}) + \beta V^\star\big(  \P^0_{F(Y^{i},\eps_1,\eps_1^0)}   \big) \Big] + \frac{\epsilon}{3}, 
\end{align}
%where $Y^{i}$ $:=$ $(\xi^i,\mra_\epsilon^i(\xi^i,U),\Lc(\xi^i,\mra_\epsilon^i(\xi^i,U))$, with $(\xi^i,U)$ $\sim$ $x^i\otimes\Uc([0,1])$.  
%by definition of $V^\star(x^i)$,  there exists $\bpi^i_\epsilon$ $\in$ $\bA$, such that    
%\begin{align} \label{Vmui} 
%V^\star(x^i) & \leq  \tilde f(x^i,\bpi_\epsilon^i) + \beta \E \big[ V^\star\big( \tilde F(x^i,\bpi_\epsilon^i,\eps_1^0)\big) \big] + \frac{\epsilon}{3}, \quad i=1,\ldots,N_\eta. 
%\end{align}
From the partition $C^\eta_k$, $k$ $=$ $1,\ldots,N_\eta$ of $\Xc^N$, associated to $\Mc^\eta$,  
we  construct the function $\bpi^\epsilon$ $:$ $\Xc^N$ $\rightarrow A^N$ as follows: we define, for all $\bx\in\Xc^N$,
\beqs
\pi^\epsilon(\bx) &=& \ba^\epsilon_k, \;\; \mbox{ when } \bx \in  C^\eta_k, \; k =1,\ldots,N^\eta.  
\enqs
Such function  $\bpi^\epsilon$  is clearly measurable. 
Let us now check that such $\bpi^\epsilon$ yields an $\epsilon$-optimal feedback policy for $\eta$ small enough.  
For $\bx$ $\in$ $\Xc^N$, we define $\bx^\eta$ $=$ $\bx_k$,  when $\bx$ $\in$ $C^\eta_k$, $k=1,\ldots,N_\eta$. Observe that $\bd_N(\bx,\bx^\eta)$ $\leq$ $\eta$. 
We then write for any $\bx$ $\in$ $\Xc^N$, 
 \begin{align}
[\Tc^{\bpi^\epsilon}_N V^\star_N](\bx) - V^\star_N(\bx) &= \; 
\Big( [\Tc^{\bpi^\epsilon}_N  V^\star_N](\bx)  - [\Tc^{\bpi^\epsilon}_N V^\star_N](\bx^\eta) \Big) + \Big([\Tc^{\bpi^\epsilon}_N V^\star_N](\bx^\eta) - V^\star_N(\bx^\eta)\Big)  \\
&  \quad\quad    + \big(V^\star_N(\bx^\eta) -  V^\star_N(\bx)\big)  \\
& \geq  \Big( [\Tc^{\bpi^\epsilon}_N V^\star_N](\bx)-[\Tc^{\bpi^\epsilon}_N V^\star_N](\bx^\eta) \Big) - \frac{\epsilon}{3} -  \frac{\epsilon}{3},  \label{interTC-N} 
\end{align}
where we used \eqref{VepsYi-N} and the fact that  $|V^\star_N(\bx^\eta) -  V^\star_N(\bx)|$ $\leq$ $\epsilon/3$ for $\eta$ small enough by uniform continuity of $V^\star_N$ in Proposition \ref{proregul-N}. 
Moreover, by observing that $\bpi^\epsilon(\bx)$ $=$ $\bpi^\epsilon(\bx^\eta)$ $=:$ $\ba$, 
we have 
\begin{align}
[\Tc^{\bpi^\epsilon}_N V^\star_N](\bx) &= \; \E \Big[ \bff(\bx,\ba) + \beta V^\star_N(\bF(\bx,\ba, \beps_1)) \Big], \\
[\Tc^{\bpi^\epsilon}_N V^\star_N](\bx^\eta) &= \; \E \Big[ \bff(\bx^\eta,\ba) + \beta V^\star_N(\bF(\bx^\eta,\ba, \beps_1)) \Big]. 
\end{align}
Under $({\bf HF_{lip}})$-$({\bf Hf_{lip}})$, and  by using  the  $\gamma$-H\"older property of $V^\star_N$ with constant $K_\star$ in Proposition \ref{proregul-N}, we then get 
\beqs
&&\big \vert [\Tc^{\bpi^\epsilon}_N V^\star_N](\bx)- [\Tc^{\bpi_\epsilon}_N V^\star_N](\bx^\eta)\big\vert \\&\leq & 2K_f \bd_N(\bx,\bx^\eta)   + \beta K_\star \E \Big[ 
\E\big[ \bd_N\big(\bF(\bx,\ba,(\eps^i_1)_{i\in\llbracket 1, N\rrbracket},e), \bF(x_\eta,\ba,(\eps^i_1)_{i\in\llbracket 1, N\rrbracket},e)\big)^\gamma\big]_{e:=\eps^0_1} \Big]  \\
&\leq & 2K_f \bd_N(\bx,\bx^\eta)   + \beta K_\star \E \Big[ 
\E\big[ \bd_N\big(\bF(\bx,\ba,(\eps^i_1)_{i\in\llbracket 1, N\rrbracket},e), \bF(\bx^\eta,\ba,(\eps^i_1)_{i\in\llbracket 1, N\rrbracket},e)\big)\big]_{e:=\eps^0_1} \Big]^\gamma  \\
&\leq & C \bd_N(\bx,\bx^\eta)^\gamma \; \leq \;  C\eta^\gamma.  
\enqs
for some constant $C$. Therefore,  $\big\vert [\Tc^{\bpi^\epsilon}_N V^\star_N](\bx)-[\Tc^{\bpi^\epsilon}_N V^\star_N](\bx ^\eta)\big\vert$ $\leq$ $\epsilon/3$, and, plugging into \eqref{interTC-N},  
we obtain $\Tc^{\bpi^\epsilon}_N V^\star_N(\bx) - V^\star_N(\bx)$ $\geq$ $-\epsilon$, for all $\bx\in\Xc^N$, which means  that $\bpi^\epsilon$ is $\epsilon$-optimal for $V^\star_N$. 
The rest of the assertions in the Theorem follows  from the verification result in Lemma  \ref{lemverif-N}.  
\ep

\bibliographystyle{plain}
\bibliography{bibliography}
 
\end{document}